\documentclass[11pt]{article} 

\usepackage{amsthm,amsfonts,amssymb,amsmath,color,mathtools}
\numberwithin{equation}{section} 
\usepackage{layout}

\usepackage{hyperref} 

\renewcommand\d{\partial}
\renewcommand\a{\alpha}

\newcommand\s{\sigma}

\newcommand\R{\mathbb R}
\newcommand\N{\mathbb N}

\def\th{\theta}

\def\l{\lambda}

\def\epsilon{\varepsilon}
\def\e{\varepsilon}

\newcommand\br{\begin{rem}}
\newcommand\er{\end{rem}}
\newcommand\bp{\begin{pmatrix}}
\newcommand\ep{\end{pmatrix}}
\newcommand\be{\begin{equation}}
\newcommand\ee{\end{equation}}
\newcommand\ba{\begin{equation}\begin{aligned}}
\newcommand\ea{\end{aligned}\end{equation}}

\newcommand\nn{\nonumber} 

\newcommand{\A}{{\mathcal A}}

\newtheorem{defi}{Definition}[section]
\newtheorem{theorem}[defi]{Theorem}
\newtheorem{proposition}[defi]{Proposition}
\newtheorem{lemma}[defi]{Lemma}
\newtheorem{corollary}[defi]{Corollary}
\newtheorem{remark}[defi]{Remark}

\newcommand{\f}[1]{\frac{1}{#1}} 
\newcommand{\Line}[2]{\underline{\hbox to #1 cm{#2}}}

\newcommand{\pa}{\partial}

\newcommand{\xwl}{\overline}

\newcommand{\Nn}{\mathbb{N}} 
\newcommand{\Zz}{\mathbb{Z}}

\renewcommand{\Re}{\mathop{\mathrm{Re}}}

\newcommand{\la}[1]{\langle #1\rangle}

\DeclareMathOperator{\xdiv}{div}

\newcommand{\norm}[1]{\lVert #1 \rVert}

\newcommand{\D}[1]{\,\mathrm{d}#1}

\begin{document}

\title{Optimal convergence rates of the  Klein-Gordon-Schr\"{o}dinger system in the nonrelativistic limit}

\author{Weizhu Bao\footnote{Department of Mathematics, National University of Singapore, 119076 Singapore; matbaowz@nus.edu.sg.}  \and Yong Lu \footnote{School of Mathematics, Nanjing University, 210093 Nanjing, China;  luyong@nju.edu.cn.}\and
    Zhiwei Zheng \footnote{School of Mathematical Sciences, Zhejiang Normal University, 321004 Jinhua, China; zhengzhiwei@zjnu.edu.cn.}
}

\date{}

\maketitle

\begin{abstract}
    In this paper, we study the Klein-Gordon-Schr\"{o}dinger system in the nonrelativistic regime $\e \to 0$, where $\e$ is proportional to the inverse of the speed of light. For initial data with limited regularity, we show that the Klein-Gordon-Schr\"{o}dinger system converges to a system of decoupled linear Schr\"odinger equations over a long time interval of order $\e^{-1}$ with error estimates of the form $(1+t)\e^2$. For more regular initial data,  we show the convergence and error estimates hold for a longer time interval of order $\e^{-\frac{4}{3}}$, where the error estimate for the Schr\"odinger component can be improved to $(1 + \e t + \e^2 t^2)\e^2$.  The specific forms of the error estimates coincide with the numerical results in \cite{BZ17}, and the $O(\e^{2})$ convergence rates coincide with the order of initial error, and thus are optimal.
\end{abstract}

\medskip

\noindent\textbf{Keywords}:  Klein-Gordon-Schr\"{o}dinger;  nonrelativistic limit;  convergence rate; WKB expansion; stability analysis in geometrical optics.




\renewcommand{\refname}{References}


\section{Introduction}

In this paper, we consider the Klein-Gordon-Schr\"{o}dinger system:
\begin{equation} \label{4-kgs}
    \begin{dcases}
        \e^2 \pa_t^2 \phi(t,x) - \Delta \phi(t,x) +\frac{\mu^2}{\e^2} \phi(t,x) - \lambda \lvert \psi(t,x) \rvert^2 =0, \\
        i \pa_t \psi(t,x) + \Delta \psi(t,x) + \lambda \phi(t,x) \psi(t,x) = 0,
    \end{dcases}
\end{equation}
where $0<\e \ll 1$ is a small parameter proportional to the inverse of the speed of light,  $\mu>0$ is the mass ratio and  $\l$ is  the coupling constant, both of which are independent of $\e$. This system generalizes the classical model of Yukawa's interaction with $\e=1$ (see Fukuda and Tsutsumi \cite{FT75,FT78}, Yukawa \cite{Y55}).
Here, the nucleons are described by the complex scalar function $\psi(t,x)$ and the mesons by the real scalar function $\phi(t,x)$, defined for $(t,x) \in \R_+ \times \R^d$. Yukawa interaction can be used to describe the strong nuclear force between nucleons, mediated by pions.

To study the dynamics of the system \eqref{4-kgs}, the initial data are usually given in a form consistent with the nonrelativistic regime (see \cite{MN02, MN03, MN08, MNO02})
\begin{equation} \label{init-data}
    \phi(0,x) = \phi_0(x), \quad \pa_t \phi(0,x) = \f{\e^2} \phi_1(x), \quad \psi(0,x) = \psi_0(x).
\end{equation}
Here, the complex-valued function $\psi_0(x)$ and the real-valued functions $\phi_0(x)$ and $\phi_1(x)$ are all assumed to be  independent of $\e$.

By multiplying ${\pa_t \phi}(t)$ and $\pa_t \psi (t)$ to \eqref{4-kgs}, we find that the Klein-Gordon-Schr\"odinger system conserves the total energy
\begin{equation} \label{unbounded-energy}
    \begin{split}
         & \int_{\R^d} \left[\f{2} \big( \e^2 |\pa_t \phi|^2  + \frac{\mu^2}{\e^2} |\phi|^2 \big)(t) + \f{2}|\nabla \phi|^2(t)+ |\nabla \psi|^2(t) - \lambda |\psi|^2 \phi (t)\right] \D{x}                 \\
         & \quad = \f{2\e^2} \int_{\R^d} \left(  |\phi_1|^2  + \mu^2 |\phi_0|^2  \right) \D{x} +  \int_{\R^d} \left(  \f{2} |\nabla \phi_0|^2+|\nabla \psi_0|^2 - \lambda |\psi_0|^2 \phi_0  \right) \D{x}.
    \end{split}
\end{equation}

For the Klein-Gordon-Schr\"odinger system \eqref{4-kgs} with parameter $\e=1$, there are extensive analytical results in the literatures \cite{BC77, BFV15, GM95, FT75, FT78, Y55,  MZ15,P04}: Fukuda and Tsutsumi \cite{FT75, FT78} showed the existence of strong solutions  in bounded domains of $\R^3$ (with smooth boundary) by means of Galerkin's method;  Baillon and  Chadam \cite{BC77} established global well-posedness and derive conservation laws in the $L^2$ framework;  Guo and Miao \cite{GM95} further studied  the global existence and the long time asymptotic behavior of the solutions.

Unlike the case where $\e=1$,  the system \eqref{4-kgs} is mathematically and numerically delicate in the nonrelativistic  regime $\varepsilon \to 0$, where the highly oscillations in time will arise and the corresponding total energy functional is unbounded and of size $ O (\e^{-2})$ (see  \eqref{unbounded-energy}).
Such nonrelativistic limit problems have attracted extensive attention in a wide range of physical equations: the Maxwell-Klein-Gordon system \cite{BMS04, MN03}, the Klein-Gordon-Zakharov system \cite{LZ25, MN08,MN10}, the Hartree equation \cite{CO06} and the Vlasov-Maxwell equation \cite{HNR18}.

The nonrelativistic limit of the nonlinear Klein-Gordon equation was systematically studied. Machihara \cite{M01} and Tsutsumi \cite{T84} showed that the nonlinear Klein-Gordon solution converges to the nonlinear Schr\"odinger solution in a finite time interval without giving any error estimates. Masmoudi and Nakanishi \cite{MN02} gave the convergence of the nonlinear Klein-Gordon equation to the Schr\"odinger equation in the energy space $H^1(\R^d)$ with certain error estimates.

Later on, the study was achieved over long time interval. Pasquali \cite{P18, P19} obtained the nonrelativistic limit valid on long time scales. In \cite{P18}, the author showed the stability of the nonrelativistic limit in the one-dimensional case for the Klein-Gordon equation with a convolution potential. In \cite{P19}, an approximation by higher-order nonlinear Schr\"odinger equations over $O(\e^{-2(r-1)})$ time was established ($r$ is the approximation order).
For quadratic and cubic defocusing Klein-Gordon equations, \cite{LZ16,BLZ24} obtained $O(\e)$ convergence over $O(\e^{-1})$ time. Moreover, for the defocusing cubic Klein-Gordon equations,  the specific form of error estimates is $(1+t)\e^{2}$ which coincides with the numerical study in \cite{BCZ14,BL26,BZ19,SZ20},  provided the initial data satisfy certain dispersive decay conditions ($H^4 \cap L^1(\mathbb{R}^3)$ as clarified in \cite{FSZ24}). Recently, along with others, Lei and Wu \cite{LW25} showed the convergence rates under weaker regularity requirements:  for initial data in $H^\alpha$ ($2 \leq \alpha \leq 4$), an error bound of  the form $\e^{2} + (\e^{2} t)^{\frac{\a}{4}}$ was given.

\medskip

Formally, following the discussion on the nonrelativistic limit of the nonlinear Klein-Gordon equation in the literatures \cite{BZ17, MN02, MNO02},  we assume the following ansatz
\begin{equation} \label{4-anstz}
    \begin{aligned}
        \mu \phi(t,x) & = e^{\frac{i\mu t}{\e^2}} \phi_{\ell}(t,x) + e^{-\frac{i\mu t}{\e^2}} \overline{\phi_{\ell}}(t,x) +  O (\e), \\
        \psi(t,x)     & = \psi_{\ell}(t,x) +  O (\e),
    \end{aligned}
\end{equation}
where $\overline{f}$ represents the complex conjugate of the function $f$.
To ensure a consistent leading-order balance, we require the $ O (1)$ oscillatory part of the residuals to vanish, leading to the following system for the profiles
\begin{equation} \label{4-sch}
    2i \partial_t \phi_{\ell} - \frac{1}{\mu} \Delta \phi_{{\ell}} = 0, \quad \partial_t \psi_{\ell} - i \Delta \psi_{{\ell}} = 0.
\end{equation}
A more rigorous justification involving higher-order correctors $(\phi_c, \psi_c)$
will be detailed in Section \ref{sec:profiles}. Corresponding to the initial data \eqref{init-data} of the original system, the equations in \eqref{4-sch} are supplemented with the following initial conditions:
\begin{equation} \label{4-sch-init}
    \phi_{\ell} (0,x) = \f{2} \big( \mu \phi_0(x) - i \phi_1(x) \big),   \quad \psi_{\ell}(0,x) = \psi_0(x).
\end{equation}

Numerical simulations were carried out in \cite{BZ17,FL26}  for \eqref{4-kgs}-\eqref{init-data} in the nonrelativistic regime $\e \to 0$ and it was observed that the exact solution is indeed well approximated by the ansatz \eqref{4-anstz} with leading profiles satisfying \eqref{4-sch}. Their computations indicate the following fact:
for initial data $(\phi_0, \phi_1, \psi_0) \in (H^3(\R^d))^3$, there holds for any $0\leq t \leq T$,
\ba\label{bz-1}
\norm{ \mu \phi (t) - \big( e^{\frac{i \mu t}{\e^2}} \phi_{\ell}(t) + e^{-\frac{i \mu t}{\e^2}} \overline{\phi_{\ell}}(t) \big) }_{H^1} + \norm{ \psi(t) - \psi_{\ell}(t) }_{H^1} \leq (C_1 + C_2 T)\e^2.
\ea
It is also observed in \cite{BZ17}  that the error estimates for the component $\psi$ achieve a uniform $\e^2$ rate:
\ba\label{bz-2}
\norm{ \psi(t) - \psi_{\ell}(t) }_{H^1} &\leq C_0 \e^2, \quad \forall \, t\in [0,T].
\ea
Our results in this paper provide a rigorous mathematical justification for these observations.

\subsection{Some useful lemmas}

As pointed out earlier in this section, formal asymptotic analysis shows that the Klein-Gordon-Schr\"odinger system converges to a system of two decoupled Schr\"odinger equations \eqref{4-sch} for the profiles $\phi_{\ell}$ and $\psi_{\ell}$. Due to the classical properties of the Schr\"odinger operator (see e.g. \cite{B99}):
\begin{equation*} 
    \norm{e^{\pm i t\Delta} f_0 }_{L^2(\R^d)}  =  \norm{ f_0 }_{L^2(\R^d)}, \quad     \norm{e^{\pm i t\Delta} f_0 }_{L^\infty(\R^d)} \leq C t^{-\frac{d}{2}} \norm{ f_0 }_{L^1(\R^d)}, \quad \forall \, t >0,
\end{equation*}
we have the following estimates:
\begin{lemma}[\cite{B99}] \label{prop:unified-est}
    Let $\s \ge 0$ and $\alpha \in \N,\, \beta \in \N^{d}$. Suppose $\phi_{\ell}$ and $\psi_{\ell}$ are the solutions to the Schr\"odinger equations in \eqref{4-sch}.  Let $f$ denote either $\phi_{\ell}$ or $\psi_{\ell}$ , and let $f_0$ be the corresponding initial datum. Then, for any $t > 0$, the following estimates hold:
    \begin{align}
        \norm{ \partial_t^{\alpha} \partial_x^{\beta} f (t) }_{H^\s(\R^d)}          & \leq C  \norm{ f_0 }_{H^{\s+2\alpha +|\beta|}(\R^d)},  \label{eq:sobolev-unified}                       \\
        \norm{ \partial_t^{\alpha} \partial_x^{\beta} f (t) }_{W^{\s,\infty}(\R^d)} & \leq C t^{-\frac{d}{2}} \norm{ f_0 }_{W^{\s+2\alpha +|\beta|, 1}(\R^d)}. \label{est-dispersive-unified}
    \end{align}
\end{lemma}

We shall also need the Kato-Ponce inequality:
\begin{lemma}[\cite{GO14, KP88}]\label{fract}
    Let $\s\ge 0$ and $1<p< \infty$ and $ 1< p_i \leq \infty$, $i =1,2,3,4$, satisfy  $\f{p} = \f{p_1} + \f{p_2} = \f{p_3} + \f{p_4}$. Then, we have
    \begin{equation*}
        \norm{ \la{\nabla}^\s (fg) }_{L^p} \leq C \left( \norm{ \la{\nabla}^\s f }_{L^{p_1}} \norm{  g }_{L^{p_2}} +  \norm{ \la{\nabla}^\s g }_{L^{p_3}} \norm{  f }_{L^{p_4}}   \right),
    \end{equation*}
    where the constant $C$ depends on $\s$, $p$ and $p_i$, $(i=1,2,3,4)$, and the operator $\la{\nabla}$ is defined by
    \begin{equation*}
        \la{\nabla} := (1 - \Delta)^{\frac{1}{2}}.
    \end{equation*}
\end{lemma}
Throughout the stability analysis, we will repeatedly use the following $H^\sigma$-bilinear estimate, which is a direct consequence of Lemma \ref{fract} with $p=2$
\begin{equation} \label{bilinear-est-app}
    \norm{ fg }_{H^\s} \leq C \left( \norm{ \la{\nabla}^\s f }_{L^{\infty}} \norm{  g }_{L^{2}} +  \norm{ \la{\nabla}^\s g }_{L^{2}} \norm{  f }_{L^{\infty}}   \right) \leq C \norm{ f}_{W^{\s, \infty}} \norm{ g }_{H^\s}.
\end{equation}
When dealing with product terms, this allows us to use the $W^{\s, \infty}$ norm to capture the $ O (t^{-\frac{d}{2}})$ decay from the  Schr\"odinger operator.

In this paper, we use $C$ to denote a generic constant that may depend on the initial data and parameters $\l$ and $\mu$, but is independent of $\e$ and $t$. While the value of $C$ may differ from line to line. We sometimes use $u(t)$ denote $u(t,\cdot)$.

\subsection{Main results}

The formal expansion presented in \eqref{4-anstz} suggests that the solutions to the equation \eqref{4-kgs} can be approximated by the solutions of the Schr\"odinger equation \eqref{4-sch}.  Our main result establishes the $O(\e^2)$ convergence rate by exploiting the dispersive decay in $d \ge 3$.
The first main result is stated as follows.

\begin{theorem} \label{main-thm-1}
    Let $d \ge 3$ and $\sigma > \frac{d}{2}$. Let the initial data satisfy
    \begin{equation} \label{init-norm-2}
        (\phi_0, \phi_1)  \in \big( H^{\s+4}(\R^d) \cap W^{\sigma, 1}(\R^d) \big)^2, \quad \psi_0 \in H^{\s+2}(\R^d) \cap W^{\sigma, 1}(\R^d).
    \end{equation}
    Let  $(\phi_\ell, \psi_\ell)$ be the solution to the Schr\"odinger system \eqref{4-sch} with initial datum \eqref{4-sch-init}.   Then, there exist $T_1 > 0$ independent of $\e$, such that the Cauchy problem \eqref{4-kgs}-\eqref{init-data} admits a unique solution $(\phi, \psi) \in C([0, \frac{T_{1}}{\e}]; {H^\s(\R^d)})$. Moreover,  the exact solution $(\phi, \psi)$ is well-approximated by the leading profile $(\phi_{\ell}, \psi_{\ell})$ with the following error estimates for any  $0\leq t\leq \frac{T_{1}}{\e}$:
    \ba\label{rate-1}
    \norm{ \mu \phi (t) - ( e^{\frac{i \mu t}{\e^2}} \phi_{\ell}(t)  +  e^{- \frac{i \mu t}{\e^2}} \overline{\phi_{\ell}}(t) ) }_{H^\s(\R^d)} +    \norm{ \psi(t) - \psi_{\ell}(t) }_{H^\s(\R^d)}  \leq  C (1+t) \e^2.
    \ea
    The constants $C$ and $T_1$ depend on the initial norms in \eqref{init-norm-2}.
\end{theorem}

The estimates in Theorem \ref{main-thm-1} establish that on the time scale $O(\e^{-1})$, both components $(\phi, \psi)$ converge at the rate of $O((1+t)\e^2)$.
However, numerical results in \cite{BZ17} suggest that for the Schr\"odinger component $\psi$, the approximation properties can be further enhanced under stronger regularity assumptions.  We can rigorously justify this observation:
\begin{theorem} \label{main-thm-2}
    In addition to the assumptions in Theorem \ref{main-thm-1}, we impose higher regularity for the initial data \eqref{init-data}:    \begin{equation} \label{init-norm-4}
        (\phi_0, \, \phi_1) \in \big( H^{\s+8}(\R^d) \cap W^{\s+4,1}(\R^d) \big)^2, \quad \psi_0 \in H^{\s+4}(\R^d) \cap W^{\s,1}(\R^d).
    \end{equation}
    Then, there exists a constant $T_2 > 0$ independent of $\e$, such that the Cauchy problem \eqref{4-kgs}-\eqref{init-data} admits a unique solution $(\phi, \psi) \in C([0, \frac{T_{2}}{\e^{\frac 43}}]; {H^\s(\R^d)})$. Moreover, there holds the following error estimates for any $0\leq t\leq \frac{T_{2}}{\e^{\frac 43}}$:
    \ba\label{rate-2}
    \norm{ \mu \phi (t) - (  e^{\frac{i \mu t}{\e^2}} \phi_{\ell}(t)  +  e^{- \frac{i \mu t}{\e^2}} \overline{\phi_{\ell}}(t) ) }_{H^\s(\R^d)} &\leq  C  (1+t + \e^2 t^2) \e^2,  \\
    \norm{ \psi(t) - \psi_{\ell}(t) }_{H^\s(\R^d)} &\leq  C (1+\e t + \e^2 t^2) \e^2 .
    \ea
    The constants $C$ and $T_2$ depend on the initial norms in \eqref{init-norm-4}.
\end{theorem}

A direct corollary is the following:
\begin{corollary}  
    Under the same assumptions as in Theorem \ref{main-thm-2}, if the time scale is restricted to $0 \leq t \leq \frac{T_{2}}{\e}$, the approximation error satisfies
    \begin{align*}
        \norm{ \mu \phi (t) - (  e^{\frac{i \mu t}{\e^2}} \phi_{\ell}(t)  +  e^{- \frac{i \mu t}{\e^2}} \overline{\phi_{\ell}}(t) ) }_{H^\s(\R^d)} & \leq  C  (1+t) \e^2, \\
        \norm{ \psi(t) - \psi_{\ell}(t) }_{H^\sigma(\R^d)}                                                                                         & \leq C \e^2.
    \end{align*}
\end{corollary}
This reveals that, over the long time scale of order $\e^{-1}$,  the error estimate for component $  \phi$ grows linearly in time as $O((1+t)\e^2)$, while the Schr\"odinger component $\psi$ achieves a uniform $ O(\e^2)$ error estimate. This result provides a rigorous justification for the numerical results obtained in \cite{BZ17} (see \eqref{bz-1} and \eqref{bz-2}).

\medskip

The rest of the paper is devoted to the proof of Theorems \ref{main-thm-1} and \ref{main-thm-2}. We will reformulate the nonrelativistic limit problems as the study of of highly oscillating solutions in nonlinear geometrical optics, and the convergence from the original Klein-Gordon-Schr\"odinger equations to the limit system of decoupled free Schr\"odinger equations corresponds to the stability of WKB approximate solutions.

\subsection{Reformulation} \label{sec:reform}

We first rewrite the equation \eqref{4-kgs} as a symmetric hyperbolic system. Following \cite{LZ16,LZ17},  we introduce a new variable
\begin{equation*} 
    \Phi := (\e (\nabla \phi)^{\rm T}, \, \e^2 \pa_t \phi,\, \mu \phi)^{\rm T}= (\e (\pa_1 \phi ,\cdots, \pa_d \phi),\, \e^2 \pa_t \phi , \, \mu \phi)^{\rm T}.
\end{equation*}
The components of $\Phi(t,x)$ are denoted by
\begin{equation*}
    \Phi(t,x) = (\Phi^{(1)}(t,x), \Phi^{(2)}(t,x), \Phi^{(3)}(t,x))^{\rm T} \in \R^{d+2},
\end{equation*}
with $\Phi^{(1)}$ being $d$-dimensional, $\Phi^{(2)}$ and $ \Phi^{(3)}$ being scalar components.

By rewriting \eqref{4-kgs} in terms of $\Phi$,  the system can be expressed as
\begin{equation} \label{4-sys-kgs}
    \begin{dcases}
        \pa_t \Phi + \f{\e} A(\pa_x) \Phi + \f{\e^2} A_0 \Phi = F(\psi, \psi), \\
        \pa_t \psi - i \Delta \psi = g(\Phi, \psi),                            \\
        \Phi(0) = \Phi(0,x), \quad  \psi(0) = \psi(0,x),
    \end{dcases}
\end{equation}
where
\begin{equation*}
    A(\pa_x) := \begin{pmatrix}
        0_{d\times d}    & -\nabla & 0_{d\times 1} \\
        -\nabla^{\rm T}  & 0       & 0             \\
        0_{1\times d}^{} & 0       & 0
    \end{pmatrix}, \quad A_0 := \begin{pmatrix}
        0_{d\times d}  & 0_{d\times 1} & 0_{d\times 1} \\
        0_{1\times d}  & 0             & \mu           \\
        0_{1 \times d} & -\mu          & 0
    \end{pmatrix},
\end{equation*}
with $0_{d\times 1}$ and $0_{d\times d}$ denoting the zero vector in $\mathbb{R}^d$ and the $d\times d$ zero matrix,
respectively, and $0_{1\times d}=(0_{d\times 1})^{\rm T}$.
We often omit these subscripts and simply write $0$ if there is no confusion in the context.
For any $\xi, \eta \in \mathbb{C}$ and $\zeta \in \mathbb{R}^{d+2}$, the bilinear terms are defined as
\begin{equation} \label{4-defF}
    F(\xi, \eta) := \begin{pmatrix} 0_{d\times 1} \\ \lambda \Re(\xi \overline{\eta}) \\ 0 \end{pmatrix}
    \quad \text{and} \quad
    g(\zeta, \eta) := \frac{i \lambda}{\mu} \zeta^{(3)} \eta,
\end{equation}
where ${\rm Re}(f)$ represents the real part of $f$.
Based on the initial data in \eqref{init-data}, the initial condition for $\Phi$ is given by
\begin{equation} \label{init-Phi0}
    \Phi(0,x) =  (\e \pa_1 \phi_0, \e \pa_2 \phi_0,\cdots, \e \pa_d \phi_0 ,\, \phi_1,\, \mu \phi_0  )^{\rm T}, \quad \psi(0,x) = \psi_0.
\end{equation}

The hyperbolic system in \eqref{4-sys-kgs} is symmetric semilinear. Due to the anti-symmetry of $A(\d_{x})$ and $A_{0}$, the large prefactor $\frac{1}{\e}$ appearing on the left-hand side of \eqref{4-sys-kgs}  has no influence on the energy estimates. Hence, with initial data in $H^\s, \s > \frac{d}{2}$, the classical existence time for the Cauchy problem \eqref{4-sys-kgs}-\eqref{init-Phi0} is $O(1)$ (see \cite{M84,M08}).

We will show the system  \eqref{4-sys-kgs} admits a global-in-time approximate solution by employing WKB expansion. The WKB expansion is a constructive method for exploring high-frequency oscillatory solutions (see e.g. \cite{R12} or \cite{JMR00}). We look for a solution $ (\Phi^a, \psi^a)$ of the form
\begin{equation} \label{4-wkb-1}
    \begin{split}
         & \Phi^a(t,x) = \sum_{n=0}^N \e^n \Phi_n(t,x), \quad \Phi_n(t,x) = \sum_{p\in \Zz} e^{i p\th} \Phi_{n,p}(t,x), \\
         & \psi^a(t,x) = \sum_{n=0}^N \e^n \psi_n(t,x),\quad \psi_n(t,x) = \sum_{p\in \Zz} e^{i p\th} \psi_{n,p}(t,x),
    \end{split}
\end{equation}
where the phase is defined as $\theta := \frac{\mu t}{\e^2}$.

It is crucial to construct the components $\Phi_{n,p}$ and $\psi_{n,p}$. This will be done step by step in Section \ref{sec:wkb}. We summarize the properties of the WKB approximate solutions in the following propositions.

\begin{proposition} \label{prop:wkb-N2}
    Under the assumptions as in Theorem \ref{main-thm-1}, there exists a WKB approximate solution $(\Phi^a, \psi^a)$ of the form \eqref{4-wkb-1} with $N=2$, solving the approximate system
    \begin{equation} \label{app-sys}
        \begin{dcases}
            \partial_t \Phi^a + \frac{1}{\e} A(\partial_x) \Phi^a + \frac{1}{\e^2} A_0 \Phi^a = F(\psi^a, \psi^a) + R_{\Phi}, \\
            \partial_t \psi^a - i\Delta \psi^a = g(\Phi^a, \psi^a) + R_{\psi},                                                \\
            \Phi^a(0,x) = \Phi(0,x) - r_{\Phi}, \quad \psi^a(0,x) = \psi(0,x) - r_{\psi},
        \end{dcases}
    \end{equation}
    for all $t>0$,   and the residual $(R_{\Phi}, R_{\psi})$ and initial error $(r_{\Phi}, r_{\psi})$ satisfy
    \begin{equation*} 
        \norm{ (R_{\Phi}, R_{\psi})(t) }_{H^\s (\R^d)} \leq C \e^2, \quad \norm{ (r_{\Phi}, r_{\psi}) }_{H^\s(\R^d)} \leq C \e^2.
    \end{equation*}
    All the components  $(\Phi_{n,p}, \psi_{n,p})$ are determined by the leading profile $(\phi_\ell, \psi_\ell)$ which satisfies the estimates for all $t>0$:
    \ba\label{init-decay}
    \| \phi_{\ell} (t) \|_{H^{\s+4}(\R^d)} +  \|   \psi_{\ell} (t) \|_{H^{\s +2}(\R^d)}    \leq C , \quad   \| (\phi_{\ell}, \psi_{\ell})(t) \|_{W^{\s,\infty}(\R^d)}  \leq C (1+t)^{-\frac{d}{2}} .
    \ea
    In particular,  the components $(\Phi_{n,p}, \psi_{n,p})$ satisfy the following estimates for all $t>0$:
    \begin{align*}
        \| \Phi_{n,p}(t) \|_{H^{\s+4-n}(\R^d)} + \| \psi_{n,p}(t) \|_{H^{\s+2}(\R^d)} \leq C, \quad \text{for all $n=0,1,2$ and $p\in \Zz$}.
    \end{align*}

\end{proposition}

We extend the construction to the higher-order case. The leading  profile  $(\phi_\ell, \psi_\ell)$ remain  the same as in Proposition \ref{prop:wkb-N2}.
To obtain a more accurate convergence rate, we introduce the corrections
\begin{align*}
    \phi_c(t) = -\frac{it}{8\mu^3} \Delta^2 \phi_{\ell}(t), \quad\text{and} \quad  \psi_c(t) = \int_0^t e^{i (t-s) \Delta} \frac{i \lambda^2}{\mu^2}  |\psi_{\ell}|^2 \psi_{\ell} (s,x) \D{s}.
\end{align*}
The following proposition summarizes the result.

\begin{proposition} \label{prop:wkb-N4}
    Under the assumptions as in Theorem \ref{main-thm-2}, there exists a WKB approximate solution $(\Phi^a, \psi^a)$ of the form \eqref{4-wkb-1} with $N=4$, solving  the approximate system \eqref{app-sys} for all $t>0$, where the residual $(R_{\Phi}, R_{\psi})$ and the initial error $(r_{\Phi}, r_{\psi})$ satisfy
    \begin{equation*} 
        \norm{ (R_{\Phi}, R_{\psi})(t) }_{H^\s(\R^d)} \leq C \e^3 + C \e^4 t, \qquad \norm{ (r_{\Phi}, r_{\psi}) }_{H^\s(\R^d)} \leq C \e^2.
    \end{equation*}
    
    All the components $(\Phi_{n,p}, \psi_{n,p})$ are determined by the leading profile $(\phi_\ell, \psi_\ell)$ and the  corrector $(\phi_c, \psi_c)$ which satisfies the following estimates for all $t>0$:
    \begin{align*} 
         & \| \phi_{\ell} (t) \|_{H^{\s+8}(\R^d)} +  \|   \psi_{\ell} (t) \|_{H^{\s + 4}(\R^d)}    \leq C , \quad  \|  \phi_{\ell} (t) \|_{W^{\s+4,\infty}(\R^d)}  + \|  \psi_{\ell} (t) \|_{W^{\s ,\infty}(\R^d)}   \leq C (1+t)^{-\frac{d}{2}} , \\
         & \| \phi_c(t) \|_{H^{\s+4}(\R^d)} + \| \psi_c(t) \|_{H^{\s+4}(\R^d)} \leq C   (1+t) ,                                                                                                                                                    \\
         & \| \phi_c(t) \|_{W^{\s,\infty}(\R^d)} \leq C   (1+t)^{-\frac{d}{2}+1}, \quad \| \psi_c(t) \|_{H^{\s}(\R^d)} \leq C (1+t)^{-d+1} .
    \end{align*}
    
    Moreover,   many of the components $(\Phi_{n,p}, \psi_{n,p})$ vanish identically, and for the nonzero ones, there hold the estimates for all $t>0$:
    \begin{align*}
         & \| \Phi_{0,\pm 1}(t) \|_{H^{\s+8}(\R^d)}  + \|  \Phi_{1,\pm 1} ( t )  \|_{H^{\s+7}(\R^d)} + \| \Phi_{2,0}(t) \|_{H^{\s+4}(\R^d)} \leq C,            \\
         &    \| \Phi_{2, \pm 1}(t)  \|_{H^{\s+4}(\R^d)} + \| \Phi_{3,0}(t) \|_{H^{\s+3}(\R^d)} +  \| (\Phi_{4,0} , \Phi_{4, \pm 1})(t)  \|_{H^{\s+2}(\R^d)} \leq C(1+t),   \\
       &  \| (\psi_{0,0}, \psi_{2, \pm 1}, \psi_{4, \pm 2}) (t) \|_{H^{\s+4}(\R^d)} \leq C, \quad       \| (\psi_{2,0}, \psi_{4,\pm 1})(t) \|_{H^{\s+4}(\R^d)} \leq C(1+t) .
    \end{align*}

\end{proposition}
These two propositions will be proved in the coming Sections \ref{prof-prop1} and \ref{prof-prop2}. We  will actually show the WKB approximate solution $(\Phi^a, \psi^a)$ is stable in the sense that for $N=2$
\begin{align} \label{eq:u-ua-2}
    \| \Phi(t) -  \Phi^a(t) \|_{H^\s(\R^d)} +  \| \psi(t)  -   \psi^a(t) \|_{H^\s(\R^d)} \leq C(1+t)\e^2, \quad \forall\,  t\in [0, \frac{T_1}{\e}],
\end{align}
and for $N=4$
\begin{align} \label{eq:u-ua-4}
    \| \Phi(t) -  \Phi^a(t) \|_{H^\s(\R^d)} + \| \psi(t) - \psi^a(t) \|_{H^\s(\R^d)}  \leq  C(1 + \e t + \e^2 t^2) \e^2, \quad \forall\,  t\in [0, \frac{T_2}{\e^{\frac{4}{3}}}].
\end{align}
Combined with the estimates for the components of the WKB solution in Propositions \ref{prop:wkb-N2} and \ref{prop:wkb-N4},
the estimate \eqref{rate-1} follows from \eqref{eq:u-ua-2}, while \eqref{rate-2} follows from \eqref{eq:u-ua-4}.

We can actually construct a WKB approximate solution up to any order $N\in \N$ as long as the initial data are smooth, and this allows us to obtain a better estimate of order $\e^{N-1}$  for $(\Phi , \psi ) -  (\Phi^a, \psi^{a}) $.  However,  the error in
\begin{equation*} 
    \big( (\Phi^a)^{(3)}(t), \psi^a(t) \big) - \big( e^{\frac{i \mu t}{\e^2}} \phi_{\ell}(t)  +  e^{- \frac{i \mu t}{\e^2}} \overline{\phi_{\ell}}(t), \psi_{\ell}(t) \big)
\end{equation*}
is at most of order $\e^{2}$, due to the fact that the WKB expansion $(\Phi^a, \psi^a)$ contains the terms $(\Phi_2, \psi_2)$ of order $\e^2$. This means that further increasing the order of the WKB expansion will not improve the final convergence rate $ O (\e^2)$. In such a sense, the $O(\e^2)$ convergence rate is optimal.

\section{Construction of WKB solutions} \label{sec:wkb}

In this section we present a step-by-step formal construction of the WKB solution. The main task is to determine the profiles $\Phi_{n,p}$ and $\psi_{n,p}$. It will be shown that the leading profiles $\phi_{\ell}$ and $ \psi_{\ell}$ are governed by the Schr\"odinger equations and all the other profiles $\Phi_{n,p}$ and $\psi_{n,p}$ are determined by $\phi_{\ell}$ and $ \psi_{\ell}$ by some induction relations.

\subsection{Formal iterative scheme}

Before proceeding with the explicit calculations, we outline the iterative procedure used to determine the WKB components. Substituting the formal expansion \eqref{4-wkb-1} into the Klein-Gordon-Schr\"odinger system \eqref{4-sys-kgs} and collecting the terms with the same powers of $\e$  and the same Fourier modes $e^{i p \th}$, we obtain the following expressions for the residuals:
\begin{equation*}
    \begin{split}
        R_{\Phi}^{} & = \sum_{n=-2}^N \sum_{p\in \mathbb{Z}} e^{i p\theta} \e^n
        \Bigl( \partial_t \Phi_{n,p} + A(\partial_x) \Phi_{n+1,p} + (A_0 + i\mu p I) \Phi_{n+2,p}
        - \bigl[ F(\psi^{a}, \psi^{a} \bigr]_{n,p} \Bigr),                      \\
        R_{\psi}^{} & = \sum_{n=-2}^N \sum_{p\in \mathbb{Z}} e^{i p\theta} \e^n
        \Bigl( \partial_t \psi_{n,p} - i \Delta \psi_{n,p} + i\mu p \psi_{n+2,p}
        - \bigl[ g(\Phi^{{a}}, \psi^{{a}}) \bigr]_{n,p} \Bigr).
    \end{split}
\end{equation*}
Here, $\left[F(\psi^a, \psi^a) \right]_{n,p}$ denotes the coefficient of $e^{i p\th} \e^n$ in the expansion of $F(\psi^a, \psi^a)$. The same notation applies to $\left[ g(\Phi^a, \psi^a) \right]_{n,p}$.
To ensure that the error between the approximate system \eqref{app-sys} and the original system \eqref{4-sys-kgs} is sufficiently small, we require the following conditions to hold for each $n \ge -2$ up to some given $N$ and all $p \in \mathbb{Z}$:
\begin{equation}\label{4-resRU}
    \begin{split}
        [{R}_{\Phi}]_{n,p} & := \pa_t \Phi_{n,p} + A(\pa_x) \Phi_{n+1,p} +(A_0 + i\mu p I) \Phi_{n+2,p} -\left[F(\psi^a, \psi^a) \right]_{n,p}  = 0, \\
        [ R_\psi]_{n,p}    & :=  \pa_t \psi_{n,p} - i \Delta \psi_{n,p} +i \mu p \psi_{n+2,p}- \left[ g(\Phi^a, \psi^a) \right]_{n,p}  = 0.
    \end{split}
\end{equation}
For convention, we set $\Phi_{n,p} = 0$ and $\psi_{n,p} = 0$ for all $n < 0$.   By the reality of $F$ and initial data and the recursive relations \eqref{4-resRU}, we impose for any $n \in \mathbb{N}$ and $p \in \mathbb{Z}$ that $\Phi_{n,-p} = \overline{\Phi}_{n,p}.$
This symmetry is consistent with the fact that the component $\Phi$ is real-valued.  Remark that $\psi_{n,p}$ is complex-valued and we do not have such symmetry on $\psi_{n,p}$.

\subsection{Derivation of specific profiles} \label{sec:profiles}
We now derive explicit expressions for the components $\Phi_{n,p}$ and $\psi_{n,p}$ by imposing the vanishing residual conditions $[R_{\Phi}]_{n,p} = 0$ and $[R_{\psi}]_{n,p} = 0$.
To ensure the existence of a solution satisfying Propositions \ref{prop:wkb-N2} and \ref{prop:wkb-N4}, we construct the profiles step by step up to order $n=2$.

\paragraph{Order $n=-2$.}
At this order, the requirement $[R_{\Phi}]_{-2,p} = 0$ yields the following algebraic  equation
\begin{equation*}
    (A_0 + i\mu p I) \Phi_{0,p} = 0, \quad \text{for each } p \in \mathbb{N}.
\end{equation*}
It follows that $\Phi_{n,p} = \Pi_{p} \Phi_{n,p}  \in \ker(A_0 + i\mu p I)$, where $\Pi_{p}$ is the orthogonal projection onto  $\ker(A_0 + i\mu p I)$.  We find that $\ker(A_0 + i\mu p I)$ is non-trivial if and only if $p\in \{0,1 \}$. Consequently, $\Phi_{0,p} = 0$ for $p \geq 2$. For the mode $p\in \{0, 1 \}$,  we have
\begin{align*}
    \ker(A_0) = \{ (v,0,0)^{\rm T}: v\in \R^{d}\};       \quad  \ker(A_0 + i\mu  I) = {\rm span}\, \{ e_{1}\},   \quad e_{1} := (0,i,1)^{\rm T}.
\end{align*}
Thus,
\ba \label{4-u01}
&  \Phi_{0,0}  = (v_0,0,0)^{\rm T} \ \mbox{for some vector-valued function $v_0(t,x)\in \R^{d}$}, \\
&  \Phi_{0,1} = \phi_{\ell}  e_{1} \ \mbox{for some scalar-valued function $\phi_{\ell}(t,x)$}.
\ea

The condition $[R_\psi]_{-2,p} = 0$ implies that $i\mu p \psi_{0,p} = 0$ for all $p \in \mathbb{N}$. Consequently,
\begin{equation} \label{4-p00}
    \psi_{0,0}(t,x) = \psi_{\ell}(t,x), \quad \psi_{0,p}(t,x) = 0 \quad \text{for all } p \neq 0,
\end{equation}
for some scalar function $\psi_{\ell}(t,x)$. The functions $\phi_{\ell}$, $v$ and $\psi_{\ell}$ will be determined at  order $n=0$.

\paragraph{Order $n=-1$.}
By considering the terms of order $ O (\e^{-1})$, the condition $[R_{\Phi}]_{-1,p} = 0$ leads to
\begin{equation*}
    A(\pa_x) \Phi_{0,p} + (A_0 + i \mu p I) \Phi_{1,p} = 0, \quad \text{for all } p \in \mathbb{N}.
\end{equation*}
This allows us to determine the components of $\Phi_{1,p}$ in terms of $\Phi_{0,p}$. Specifically,
\begin{align}
    \Phi_{1,0} & = (v_1, 0, \f{\mu} \xdiv v_0)^{\rm{T} },  \quad  \text{for some vector-valued function $v_1 \in \R^{d}$}, \notag             \\
    \Phi_{1,1} & = (\f{\mu} (\nabla \phi_{\ell})^{\rm T}, 0,0 ) + w_1 e_1, \quad  \text{for some scalar-valued function $w_1$}, \label{4-u11} \\
    \Phi_{1,p} & = 0, \quad \text{for all $p \in \Nn$ such that }  p\ne 1. \notag
\end{align}

From  $[R_\psi]_{-1,p} = 0$, we find
\begin{equation*}
    i \mu \psi_{1,p} = 0, \quad \text{for all } p\in \Nn.
\end{equation*}
It follows that $\psi_{1,p} = 0$ for all $p \ne 0$, and then the only non-trivial component is $\psi_{1,0}$.

\paragraph{Order $n=0$: determination of $\phi_{\ell}$ and $\psi_{\ell}$.}

By requiring  $[R_\Phi^a]_{0,p} = 0$, we obtain the following equation for the profile $\Phi_{2,p}$:
\begin{equation} \label{Phi0}
    \pa_t \Phi_{0,p} + A(\pa_x) \Phi_{1,p} + (A_0 + i\mu p I) \Phi_{2,p} = \left[ F(\psi^a, \psi^a) \right]_{0,p}.
\end{equation}
Utilizing the definition of $F$ in \eqref{4-defF}, the zero-th mode of the nonlinear term is
\begin{equation*}
    \left[ F(\psi^a, \psi^a) \right]_{0,0} = (0, \lambda |\psi_{0,0}|^2,0)^{\rm{T}} = (0,\lambda |\psi_{\ell}|^2,0)^{\rm{T}}.
\end{equation*}
Consequently, the components of the second-order profile $\Phi_2$ are
\begin{equation} \label{4-u2p}
    \Phi_{2,0} = \begin{pmatrix}
        v_2 \\
        0   \\
        \frac{\lambda}{\mu} |\psi_{\ell}|^2 + \f{\mu} \xdiv v_1
    \end{pmatrix}, \quad
    \Phi_{2,1} = \begin{pmatrix}
        \f{\mu} \nabla w_1                   \\
        \frac{1}{\mu} \partial_t \phi_{\ell} \\
        0
    \end{pmatrix} + \phi_c e_1, \quad
    \Phi_{2,p} = 0, \quad \text{for all } p \ge 2,
\end{equation}
for some scalar-valued function $\phi_{c}$ and some vector-valued function $v_{2}\in \R^{d}$.

For zero mode $p=0$, by \eqref{4-u01}, \eqref{4-u11}, \eqref{Phi0} and \eqref{4-u2p}, we deduce that
\begin{align}
    \pa_t v_0 = 0. \label{eq:v0}
\end{align}
As will be seen later (see the argument below \eqref{eq:e2}), we could set initial datum $v_0(0,x) = 0$ and thus $v_0(t,x) \equiv 0$,  and consequently $\Phi_{0,0} = 0$.  By induction, we shall always exclude the zero mode and set $\Phi_{n,0} = 0$ for any $n$.

We now turn to consider the mode $p=1$, which yields the equation for the profile $\phi_\ell$.  Specifically, substituting the expressions for $\Phi_{0,1}$ and $\Phi_{1,1}$ from \eqref{4-u01} and \eqref{4-u11} into \eqref{Phi0} implies
\begin{equation*}
    \begin{dcases}
        i \mu \Phi_{2,1}^{(1)} = 0,                                                                                 \\
        i \pa_t \phi_{\ell} - \frac{1}{\mu} \Delta \phi_{\ell} + i \mu \Phi_{2,1}^{(2)} + \mu \Phi_{2,1}^{(3)} = 0, \\
        \pa_t \phi_{\ell} - \mu \Phi_{2,1}^{(2)} + i \mu \Phi_{2,1}^{(3)} = 0.
    \end{dcases}
\end{equation*}
By eliminating the terms $\Phi_{2,1}^{(2)}$ and $\Phi_{2,1}^{(3)}$, we deduce the free Schr\"odinger equation in $\phi_{\ell}$:
\begin{equation} \label{eq-for-z}
    2i \pa_t \phi_{\ell} - \frac{1}{\mu} \Delta \phi_{\ell} = 0.
\end{equation}

Regarding the equation for $\psi_\ell$, the condition $[R_\psi]_{0,p} = 0$ yields
\begin{equation} \label{4-eqp2}
    \pa_t \psi_{0,p} - i \Delta \psi_{0,p} + i\mu p \psi_{2,p} = \left[ g(\Phi^a, \psi^a) \right]_{0,p}.
\end{equation}
By substituting the form of $g$ from \eqref{4-defF} into the right-hand side, we obtain
\begin{equation}\label{g0-1}
    \left[ g(\Phi^a, \psi^a) \right]_{0} = \frac{i\lambda}{\mu} \phi_{\ell} \psi_{\ell} e^{i\theta} + \frac{i\lambda}{\mu} \overline{\phi_{\ell}} \psi_{\ell} e^{-i\theta}.
\end{equation}
Specifically, when $p=0$, the nonlinear source term vanishes since it only contains modes $p=\pm 1$. Substituting the profiles from \eqref{4-p00} into \eqref{4-eqp2}, we obtain
\begin{equation*} 
    \partial_t \psi_{\ell} - i \Delta \psi_{\ell} = 0.
\end{equation*}
For mode $p=0$, we also introduce the second-order correction $\psi_{c}$ (previously $\psi_{2,0}$), which is yet unknown.
It will be determined by the conditions arising at  step $n=2$.

For mode $p=\pm 1$, we deduce from \eqref{4-eqp2} and \eqref{g0-1} that
\ba \label{eq:psi2-1}
&  i \mu \psi_{2,1} = \frac{i\lambda}{\mu} \phi_{\ell} \psi_{\ell}, \quad \text{leading to} \quad \psi_{2,1} = \frac{\lambda}{\mu^2} \phi_{\ell} \psi_{\ell},
\\
& -i \mu \psi_{2,-1} = \frac{i\lambda}{\mu} \overline{\phi_{\ell}} \psi_{\ell}, \quad \text{leading to} \quad \psi_{2, -1} = -\frac{\lambda}{\mu^2} \overline{\phi_{\ell}} \psi_{\ell}.
\ea
These relations show that the higher-order corrections $\psi_{2,\pm 1}$ are determined algebraically by the leading-order profiles $\phi_{\ell}$ and $\psi_{\ell}$.

Summarizing the derivation from \eqref{4-u01} to \eqref{eq:psi2-1}, the approximate solutions for $\Phi$ and $\psi$ up to $ O (\e^2)$ can be expressed as follows:
\begin{align} \label{2-wkb-Phi-psi}
    \begin{split}
        \Phi^a & = \begin{pmatrix}
                       v_0                    \\
                       i \phi_{\ell} e^{i\th} \\
                       \phi_{\ell} e^{i \th}
                   \end{pmatrix}
        + \e \begin{pmatrix}
                 v_1 + \frac{1}{\mu} \nabla \phi_{\ell} e^{i\th } \\
                 i w_{1} e^{i\th}                                 \\
                 \f{\mu}\xdiv v_0 + w_{1} e^{i\th}
             \end{pmatrix}
        + \e^2  \begin{pmatrix}
                    v_2 + \f{\mu} \nabla w_1 e^{i \th}                            \\
                    \frac{1}{\mu} \pa_t \phi_{\ell} e^{i\th} + i \phi_{c}e^{i\th} \\
                    \frac{\lambda}{\mu} | \psi_{\ell}|^2 + \f{\mu} \xdiv v_1 +  \phi_{c} e^{i\th}
                \end{pmatrix}   + \text{c.c.},                                                                                                               \\
        \psi^a & = \psi_{\ell} +\e \psi_1 +\e^2 \big(\psi_{c} + \frac{\lambda}{\mu^2} \phi_{\ell} \psi_{\ell} e^{i \theta} - \frac{\lambda}{\mu^2} \overline{\phi_{\ell} }\psi_{\ell} e^{- i\theta} \big).
    \end{split}
\end{align}
where $\text{c.c.}$ denotes the complex conjugate of the terms in front of it.

The initial errors are given by
\begin{equation*} 
    r_{\Phi}(x) = \Phi(0,x) - \Phi^a(0,x), \quad r_{\psi}(x) = \psi(0,x) - \psi^a(0,x),
\end{equation*}
with $\Phi(0,x)$ and $\psi(0,x)$ as in \eqref{init-Phi0}.
Thus, to minimize the initial error, we impose
\begin{align} \label{eq:e2}
    \begin{dcases}
        v_0(0,x) + \e ( v_1(0,x) + \f{\mu} \nabla \phi_\ell(0,x) + \f{\mu} \nabla \xwl{\phi_\ell}(0,x) ) = \e \nabla \phi_0(x), \\
        i \phi_\ell(0,x) - i \xwl{\phi_\ell}(0,x) = \phi_1(x),                                                                  \\
        \phi_\ell(0,x) + \xwl{\phi_\ell}(0,x) = \mu \phi_0(x),                                                                  \\
        \psi_\ell(0,x) = \psi_0(x).
    \end{dcases}
\end{align}
It follows directly that the initial data for $\phi_{\ell}$ and $\psi_{\ell}$ are given by \eqref{4-sch-init}, while the initial data for $v_0$ and $v_1$ are determined by the first equation in \eqref{eq:e2} so that
\begin{align} \label{init:v0-v1}
    v_0(0,x)  \equiv 0, \quad  v_1(0,x)  \equiv 0.
\end{align}
Since we have $\pa_t v_0 =0$ in \eqref{eq:v0}, we have $v_0(t,x) \equiv 0$.

We set the remaining initial data to zero:
\begin{align} \label{init:zero}
    w_{1}(0,x) \equiv 0, \quad v_2(0,x) \equiv0,   \quad \phi_c(0,x) \equiv 0, \quad \psi_{1,0}(0,x)\equiv0, \quad \psi_c(0,x) \equiv 0.
\end{align}
This choice is consistent with the initial data of the original system \eqref{4-sys-kgs} and does not affect the convergence rate of the WKB approximation.

\paragraph{Order $n=1$.}
Setting  $[R_{\Phi}]_{1,p} = 0$ gives the following equation for $\Phi_{3,p}$:
\begin{align} \label{eq:Phi3p}
    \pa_t \Phi_{1,p} + A(\pa_x) \Phi_{2,p} + (A_0 + i\mu p I) \Phi_{3,p} = \left[ F(\psi^a, \psi^a) \right]_{1,p}.
\end{align}
For mode $p=0$, we have $\pa_t v_1 = 0$. Together with the initial data in \eqref{init:v0-v1}, this implies $v_1(t,x) \equiv 0$. Moreover, from \eqref{eq:Phi3p} and \eqref{4-u2p}, we obtain
\begin{align*}
    \Phi_{3,0} = \begin{pmatrix}
                     v_3 \\
                     0   \\
                     \frac{2 \l}{\mu} \Re (\psi_{1} \overline{\psi_{1} } ) + \f{\mu} \xdiv v_2
                 \end{pmatrix}, \quad \mbox{for some vector-valued function $v_{3}\in \R^{d}$}.
\end{align*}

For mode $p=1$, it follows directly from \eqref{eq:Phi3p} and \eqref{4-u2p} that
\begin{align*}
    \begin{dcases}
        -i \nabla \phi_c + i \mu \Phi_{3,1}^{(1)} = 0,                                                \\
        i \pa_t w_{1} - \frac{1}{\mu} \Delta w_1 + i \mu \Phi_{3,1}^{(2)} + \mu \Phi_{3,1}^{(3)} = 0, \\
        \pa_t w_{1} - \mu \Phi_{3,1}^{(2)} + i \mu \Phi_{3,1}^{(3)} = 0.
    \end{dcases}
\end{align*}
Thus, we have $2i \pa_t w_1 - \f{\mu} \Delta w_1 = 0$. Combining with the initial data in \eqref{init:zero}, we have $w_1(t,x) \equiv 0$. Moreover, we have the expression for $\Phi_{3,1}$:
\begin{align*}
    \Phi_{3,1} = \begin{pmatrix}
                     \frac{1}{\mu} \nabla \phi_c \\
                     0                           \\
                     0
                 \end{pmatrix} + w_3 e_1, \quad \text{for some scalar-valued function $w_3$}.
\end{align*}

For $p \ge 2$, since the right-hand side of \eqref{eq:Phi3p} vanishes and $(A_0 + i\mu p)$ is invertible,  we obtain $\Phi_{3,p} = 0$. Consequently,  we have the expression for $\Phi_3$:
\begin{align*}
    \Phi_3 = \begin{pmatrix}
                 v_3 + \frac{1}{\mu} \nabla \phi_{c} e^{i\th} \\
                 w_3 e^{i\th}                                 \\
                 \frac{2\lambda}{\mu} \Re (\psi_{1} \overline{\psi_{1}}) + \f{\mu} \xdiv v_2 +  i w_3 e^{i\th}
             \end{pmatrix} + \text{c.c.}.
\end{align*}

Equations $[R_{\psi}]_{1,p} =0$ can be written as
\ba\label{R-psi-1p}
\pa_t \psi_{1,p} - i \Delta \psi_{1,p} + i \mu p \psi_{3,p} = \left[ g(\Phi^a, \psi^a) \right]_{1,p}.
\ea
By \eqref{2-wkb-Phi-psi}, using the fact $v_0\equiv 0$, $v_{1} \equiv 0$ and $w_1=0$, we obtain
\begin{align} \label{eq:g1}
    [ g(\Phi^a, \psi^a) ]_{1} = \frac{i\l}{\mu} \phi_\ell \psi_{1,0} e^{i\th} + \frac{i \l}{\mu} \overline{\phi_\ell} \psi_{1,0} e^{-i\th}.
\end{align}
Thus, for $p=0$, we deduce from \eqref{R-psi-1p} and \eqref{eq:g1} that
$$\pa_t \psi_{1,0} - i \Delta \psi_{1,0} = 0.$$
Combining with the zero initial data in \eqref{init:zero}, we have $\psi_{1,0}(t,x) \equiv 0$.  Returning to \eqref{eq:g1}, we have $[ g(\Phi^a, \psi^a) ]_{1} = 0$. Thus,   $\psi_{3,1} = 0$ for all $p \neq 0$. Therefore, we have
\begin{align*}
    \psi_3(t,x) = \psi_{3,0}(t,x), \quad \text{for some scalar-valued function $\psi_{3,0}$}.
\end{align*}

In the computation of $\Phi_{3}$ and $\psi_{3}$, we introduced new profiles $v_3$, $w_3$ and $\psi_{3,0}$. However, as mentioned in Section \ref{sec:reform},
even including these $O(\e^3)$ order profiles in the construction of WKB solutions,  the final convergence rate $O(\e^2)$ will not be improved.  Hence, we simply set them to be zero and obtain
\begin{equation*}
    \Phi_3 = \begin{pmatrix}
        \frac{1}{\mu} \nabla \phi_{c} \\
        0                             \\
        0
    \end{pmatrix} e^{i \theta} + \text{c.c.}, \quad \psi_3 = 0.
\end{equation*}

\paragraph{Order $n=2$.} In this order, the equations for profiles $\phi_c$ and $\psi_c$ will be derived.  As in order $n=1$, we shall not introduce new unknown $O(\e^{4})$ profiles and simply impose $\Pi_{0}\Phi_{4,0} =\Phi_{4,0}^{(1)} =  0$, $\Pi_{1}\Phi_{4,1} = 0$ and $\psi_{4,0}=0$.

Based on the profiles $\psi_{2, \pm 1}$ constructed in   \eqref{eq:psi2-1}, we find
\begin{equation*}
    \begin{split}
        [F(\psi^a, \psi^a)]_{2}^{(2)} & = 2 \lambda \Re \big( \overline{\psi_{\ell}} \big( \psi_{c} + \psi_{2,1} e^{i \theta} + {\psi_{2,-1}} e^{-i \theta} \big)   \big)                                                                \\
                                      & = 2 \lambda \Re(\overline{\psi_{\ell}} \psi_{c}) + 2\lambda \Re\big(\frac{\l}{\mu^2}\phi_{\ell} |\psi_{\ell}|^2 e^{i\th} - \frac{\l}{\mu^2}\overline{\phi_{\ell}} |\psi_{\ell}|^2 e^{-i\th}\big) \\
                                      & = 2 \lambda \Re(\psi_{\ell} \overline{\psi_{c}}).
    \end{split}
\end{equation*}
By setting $[R_{\Phi}]_{2, p} = 0$, we obtain
\begin{align*}
    \pa_t \Phi_{2,p} + A(\pa_x) \Phi_{3,p} + (A_0 + i\mu p I) \Phi_{4,p} = [F(\psi^a, \psi^a)]_{2,p}.
\end{align*}
For mode $p=0$, we have
\begin{align} \label{eq:Phi40}
    \begin{cases}
        \pa_t v_2 = 0,                                                         \\
        \mu \Phi_{4,0}^{(3)} = 2 \lambda \Re(\psi_{\ell} \overline{\psi_{c}}), \\
        \frac{\l}{\mu}\pa_t |\psi_{\ell}|^2 -  \mu \Phi_{4,0}^{(2)} = 0.
    \end{cases}
\end{align}
Similarly, for mode $p=1$, we have
\begin{align} \label{eq:Phi41}
    \begin{cases}
        i \mu \Phi_{4,1}^{(1)} = 0,                                                                                                   \\
        \f{\mu} \pa_t^2 \phi_\ell + i \pa_t \phi_c - \frac{1}{\mu} \Delta \phi_c + i \mu \Phi_{4,1}^{(2)} + \mu \Phi_{4,1}^{(3)} = 0, \\
        \pa_t \phi_c - \mu \Phi_{4,1}^{(2)} + i \mu \Phi_{4,1}^{(3)} = 0.
    \end{cases}
\end{align}
From \eqref{eq:Phi40} and \eqref{eq:Phi41}, we determine that the component $\Phi_{4}$ is given by
\begin{equation*}
    \Phi_4 = \begin{pmatrix}
        0                                           \\
        \frac{\lambda}{\mu^2} \pa_t |\psi_{\ell}|^2 \\
        \frac{2\lambda}{\mu} \Re (\psi_{\ell} \overline{\psi_{c}})
    \end{pmatrix} + \begin{pmatrix}
        0                            \\
        \frac{1}{\mu} \pa_t \phi_{c} \\
        0
    \end{pmatrix} e^{ i \theta} + \text{c.c.}.
\end{equation*}
We observe that the equation for $v_2$ in \eqref{eq:Phi40} implies $v_2(t,x) \equiv 0.$
From the mode $p=1$,
a non-homogeneous Schr\"odinger equation is obtained:
\begin{equation*} 
    2 i \pa_t \phi_{c} - \frac{1}{\mu} \Delta \phi_{c} + \frac{1}{\mu} \pa_t^2 \phi_{\ell} = 0.
\end{equation*}
Together with the zero initial condition \eqref{init:zero}, we obtain
\begin{equation} \label{formula-phic}
    \phi_{c}(t,x) = \int_0^t e^{\frac{-i (t-s)}{2\mu} \Delta}  \frac{i}{2\mu} \pa_s^2 \phi_{\ell} (s,x) \D{s}.
\end{equation}
From \eqref{eq-for-z}, it follows that
\begin{align*}
    e^{\frac{-i (t-s)}{2\mu} \Delta} \pa_s^2 \phi_{\ell} (s,x)
     & = -\frac{1}{4\mu^2} e^{\frac{-i t}{2\mu} \Delta}  \Delta^2 \phi_{\ell}(0,x)   = -\frac{1}{4\mu^2} \Delta^2 \phi_{\ell}(t,x).
\end{align*}
Thus, we deduce the formula for $\phi_{c}$ from \eqref{formula-phic}:
\begin{align*} 
    \phi_{c}(t,x) = -\frac{it}{8\mu^3} \Delta^2 \phi_{\ell}(t,x).
\end{align*}
It is clear from this formula that $\phi_{c}$ is of order $O(t)$ in $H^\s(\mathbb{R}^d)$ norm.

To construct the higher-order profile $\psi_4$, we first compute the nonlinear term $[g(\Phi^a, \psi^a)]_2$. By substituting the profiles $\Phi_0$, $\psi_2$, $\Phi_2$ and $\psi_0$ into the expression $g$, we obtain
\begin{align} \label{eq:last-g}
    \begin{split}
         & [g(\Phi^a, \psi^a)]_{2}      = \frac{i \lambda}{\mu} \big[ \big( \phi_{\ell} e^{i \theta} + \overline{\phi_{\ell}} e^{-i \theta} \big) \big( \psi_{c} + \frac{\l}{\mu^2} \phi_{\ell} \psi_{\ell} e^{i \theta} - \frac{\l}{\mu^2} \overline{\phi_{\ell}} \psi_{\ell} e^{-i\theta} \big) + \big( \frac{\lambda}{\mu} |\psi_{\ell}|^2 + \phi_{c} e^{i\theta} + \overline{\phi_{c}} e^{-i \theta} \big) \psi_{\ell} \big]     \\
         & \quad  = \frac{i\lambda^2}{\mu^2} |\psi_{\ell}|^2 \psi_{\ell}  + \frac{i \l}{\mu} \left( \phi_{\ell} \psi_{c} + \phi_{c} \psi_{\ell} \right) e^{i \theta} + \frac{i\l^2}{\mu^3} \phi_{\ell}^2 \psi_{\ell} e^{2 i\theta} + \frac{i\lambda}{\mu} \left( \overline{\phi_{\ell}} \psi_{c} + \overline{\phi_{c}} \psi_{\ell} \right) e^{-i \theta} - \frac{i\l^2}{\mu^3} \overline{\phi_{\ell}}^2 \psi_{{\ell}} e^{-2i\theta}.
    \end{split}
\end{align}
The equation $(R_\psi)_{2,p}=0$ yields
\begin{align} \label{eq:psi2p}
    \pa_t \psi_{2,p} - i \Delta \psi_{2,p} + i \mu p \psi_{4,p} = [g(\Phi^a, \psi^a)]_{2,p}.
\end{align}

By \eqref{eq:last-g}, we find that the right-hand side of \eqref{eq:psi2p} only contains modes $p=0,\, \pm 1,\, \pm 2$.  For mode $p=0$, we deduce the equation in  profile $\psi_{c}$ from \eqref{eq:psi2p}:
\begin{align*}
    \pa_t \psi_{c} - i \Delta \psi_{c} = \frac{i \lambda^2}{\mu^2} |\psi_{\ell}|^2 \psi_{\ell}.
\end{align*}
Combined with the initial condition \eqref{init:zero}, this yields
\begin{equation} \label{formula-psic}
    \begin{split}
        \psi_{c}(t,x) = \int_0^t e^{i (t-s) \Delta} \frac{i \lambda^2}{\mu^2}  |\psi_{\ell}|^2 \psi_{\ell} (s,x) \D{s}.
    \end{split}
\end{equation}

As before, we will simply set $\psi_{4,0} = 0$. For the remaining  modes $p=\pm 1$ and $p=\pm 2$, the right-hand side of \eqref{eq:psi2p} is explicitly given by \eqref{eq:last-g}, and the factor $i \mu p$ is invertible. Together with the fact $\psi_{2,p}=0$ for $p\neq 0, \pm 1$,  we thus have
\begin{align*}
    \psi_4 = \psi_{4,1} e^{i \theta} + \psi_{4,-1} e^{-i \theta} + \psi_{4,2} e^{2 i \theta} + \psi_{4,-2} e^{-2i \theta},
\end{align*}
where the profiles $\psi_{4,\pm 1}$ and $\psi_{4,\pm 2}$ are determined algebraically from \eqref{eq:psi2p}:
\begin{align*}
    \psi_{4,1}  & = \frac{\lambda}{\mu^2} \left( \phi_{\ell} \psi_{c} + \phi_{c} \psi_{\ell} \right) + \frac{i \l}{\mu^3} \left( \pa_t (\phi_{\ell} \psi_{{\ell}}) - i \Delta (\phi_{\ell} \psi_{{\ell}}) \right),                                             \\
    \psi_{4,-1} & = -\frac{\lambda}{\mu^2} \left( \overline{\phi_{\ell}} \psi_{c} + \overline{\phi_{c}} \psi_{\ell} \right) - \frac{i\l}{\mu^3} \left( \pa_t (\overline{\phi_{{\ell}}} \psi_{\ell}) - i \Delta (\overline{\phi_{{\ell}}} \psi_{\ell}) \right), \\
    \psi_{4,2}  & = \frac{\l^2}{2\mu^4} \phi_{\ell}^2 \psi_{\ell}, \quad \psi_{4,-2} = \frac{\l^2}{2\mu^4} \overline{\phi_{\ell}}^2 \psi_{{\ell}}.
\end{align*}

Finally, the higher-order WKB approximate solution is given by
\begin{align} \label{wkb-Phi-final}
    \begin{split}
        \Phi^a & = \begin{pmatrix}
                       0                            \\
                       i \phi_{{\ell}} e^{i \theta} \\
                       \phi_{{\ell}} e^{i \theta}
                   \end{pmatrix} + \e \begin{pmatrix}
                                          \f{\mu} \nabla \phi_\ell e^{i\th} \\
                                          0                                 \\
                                          0
                                      \end{pmatrix} + \e^2 \begin{pmatrix}
                                                               0                                                     \\
                                                               \f{\mu} \pa_t \phi_\ell e^{i\th} + i  \phi_c e^{i\th} \\
                                                               \frac{\l}{\mu} |\psi_\ell|^2 + \phi_c e^{i\th}
                                                           \end{pmatrix}  + \e^3 \begin{pmatrix}
                                                                                     \f{\mu} \nabla \phi_c e^{i\th} \\
                                                                                     0                              \\
                                                                                     0
                                                                                 \end{pmatrix} \\
               & \quad  + \e^4 \begin{pmatrix}
                                   0                                                                    \\
                                   \f{\mu} \pa_t \phi_c e^{i\th} + \frac{\l}{\mu^2} \pa_t |\psi_\ell|^2 \\
                                   \frac{2\l}{\mu} \Re(\psi_\ell \xwl{\psi_c})
                               \end{pmatrix} + \text{c.c.},
    \end{split}
\end{align}
and
\begin{equation} \label{wkb-psi-final}
    \psi^a = \psi_{{\ell}} + \e^2 \big(\psi_{c} + \frac{\l}{\mu^2} \phi_{\ell} \psi_{\ell} e^{i \theta} -\frac{\l}{\mu^2} \overline{\phi_{\ell}} \psi_{\ell} e^{- i\theta} \big) + \e^4 \psi_4.
\end{equation}

With these constructions, the approximate solution $(\Phi^a, \psi^a)$ is well-defined, providing a foundation for the subsequent error analysis and convergence estimates.

\subsection{Low regularity case: proof of Proposition \ref{prop:wkb-N2}} \label{prof-prop1}

We present the residual estimates for the WKB approximations at order $N=2$. The construction relies on the profiles derived in Section \ref{sec:profiles}.

Up to order $N=2$, we define the approximate solution $(\Phi^a, \psi^a)$ by adding the first algebraic corrections to the leading-order profiles:
\begin{equation} \label{sche-2}
    \Phi^a = \begin{pmatrix}
        \frac{\e}{\mu} \nabla \phi_{{\ell}} e^{i \theta}                                  \\
        i \phi_{{\ell}} e^{i \theta} + \frac{\e^2}{\mu} \pa_t \phi_{{\ell}} e^{ i \theta} \\
        \phi_{{\ell}} e^{i \theta} + \frac{\e^2 \lambda}{\mu} |\psi_{{\ell}} |^2
    \end{pmatrix} + \text{c.c.}, \quad   \psi^a = \psi_{{\ell}} + \e^2 \big( \frac{\lambda}{\mu^2} \phi_{\ell} \psi_{\ell}  e^{i \theta} - \frac{\l}{\mu^2} \overline{\phi_{\ell}} \psi_{\ell} e^{- i\theta} \big).
\end{equation}
This ansatz is accurate up to $O(\e^2)$ in the Klein-Gordon-Schr\"odinger system.

The uniform estimates of the components $\phi_\ell$ and $\psi_\ell$ follow directly from \eqref{eq:sobolev-unified} and \eqref{est-dispersive-unified}. Given initial data satisfying \eqref{init-norm-2}, for all $t>0$,
we have
\begin{align*}
    \| \phi_\ell(t) \|_{H^{\s+4}}  \leq C  \| ( \phi_0, \phi_1) \|_{H^{\s+4}}, \quad \| \psi_\ell(t) \|_{H^{\s+2}} \leq C \| \psi_0 \|_{H^{\s+2} },
\end{align*}
and
\begin{align*}
    \| \phi_\ell(t) \|_{W^{\s, \infty}} \leq C (1+t)^{-\frac{d}{2}} \| (\phi_0, \phi_1) \|_{W^{\s,1}\cap H^\s}, \quad \| \psi_\ell (t) \|_{W^{\s, \infty}} \leq C (1+t)^{-\frac{d}{2}} \| \psi_0 \|_{W^{\s,1}\cap H^\s}.
\end{align*}

Since $\Phi^a$ and $\psi^a$ are differential polynomials in $\phi_\ell$ and $\psi_\ell$, the algebra property of $H^\s$ yields the following uniform bounds for all nonzero compoents:
\begin{align*}
    \begin{split}
        \| \Phi_{0,1}(t) \|_{H^{\s+4}} & \leq C \| \phi_\ell(t) \|_{H^{\s+4}} \leq C \| (\phi_0, \phi_1) \|_{H^{\s+4}},         \\
        \| \Phi_{1,1}(t) \|_{H^{\s+3}} & \leq C \|  \nabla \phi_\ell(t) \|_{H^{\s+3}} \leq C \| (\phi_0, \phi_1) \|_{H^{\s+4}}, \\
        \| \Phi_{2,0}(t) \|_{H^{\s+2}} & \leq  C \| \psi_\ell(t) \|_{H^{\s+ 2}}^2 \leq   C \| \psi_0 \|_{H^{\s+2}}^2,           \\
        \| \Phi_{2,1}(t) \|_{H^{\s+2}} & \leq C \| \pa_t \phi_\ell(t) \|_{H^{\s+2}} \leq  C \| (\phi_0, \phi_1) \|_{H^{\s+4}},
    \end{split}
\end{align*}
and  \begin{align*}
    \begin{split}
         & \| \psi_{0,0}(t) \|_{H^{\s+2}}                                  \leq  \| \psi_\ell(t) \|_{H^{\s+2}} \leq C \| \psi_0 \|_{H^{\s+2}}, \\
         & \| \psi_{2,1}(t) \|_{H^{\s+2}} + \| \psi_{2,-1}(t) \|_{H^{\s+2}} \leq C \| (\phi_0, \phi_1) \|_{H^{\s+2}} \| \psi_0 \|_{H^{\s+2}}.
    \end{split}
\end{align*}
Therefore, the WKB approximation $(\Phi^a, \psi^a)$ is uniformly bounded in $H^{\s+2}(\R^d)$ for all $t \ge 0$.
We require the profiles to be in the Sobolev space $H^{\s+2}(\R^d)$ to ensure that the residuals lie in $H^\s(\R^d)$, which is necessary for the error estimates in Theorem \ref{main-thm-1}.

From the WKB expansion in Section \ref{sec:profiles},  it follows that the approximate solution $(\Phi^{a}, \psi^{a})$ defined in \eqref{sche-2} satisfies the approximate solution \eqref{app-sys} with the Klein-Gordon residual
\begin{equation*} 
    R_{\Phi} = \e^2 \begin{pmatrix}
        0                                                 \\
        \frac{2}{\mu}   \Re(\pa_t^2 \phi_{\ell} e^{i\th}) \\
        \frac{\lambda}{\mu} \pa_t |\psi_{\ell}|^2
    \end{pmatrix} - \lambda \e^4 \begin{pmatrix}
        0          \\
        |\psi_2|^2 \\
        0
    \end{pmatrix},
\end{equation*}
and the Schr\"odinger residual
\begin{equation*}
    \begin{split}
        R_{\psi} & =  \e^2 \big( \frac{\lambda}{\mu^2} \partial_t (\phi_{\ell} \psi_{\ell}) e^{i \theta} - \frac{\lambda}{\mu^2} \partial_t (\overline{\phi_{\ell}} \psi_{\ell}) e^{-i\theta}  - i \Delta \psi_2 \big) \\
                 & \quad  - \e^2 \frac{i \lambda }{\mu} \big( \frac{\lambda}{\mu}|\psi_{\ell}|^2 \psi_{\ell}
        + \psi_2 (\phi_{\ell} e^{i\theta} + \overline{\phi_{\ell}}e^{-i \theta} - \frac{\lambda \e^2}{\mu} |\psi_{\ell}|^2 ) \big).
    \end{split}
\end{equation*}
By \eqref{eq:sobolev-unified}, we obtain that for all $t \ge 0$,
\begin{align} \label{est:Rphi-2}
    \begin{split}
        \| R_{\Phi}(t) \|_{H^{\s}} & \leq C \e^2 \| \pa_t^2 \phi_\ell \|_{H^\s} + C\e^2 \| \pa_t |\psi_\ell|^2 \|_{H^\s} + C \e^4 \| \psi_2 \|_{H^\s}^2                                           \\
                                   & \leq C \e^2\big( \| (\phi_0, \phi_1) \|_{H^{\s+4}} +   \| \psi_0 \|_{H^{\s+2}}^2 +   \e^2\| (\phi_0, \phi_1) \|_{H^{\s+4}}^2 \| \psi_0 \|_{H^{\s+2}}^2\big),
    \end{split}
\end{align} 
and
\begin{align} \label{est:Rpsi-2}
    \begin{split}
        \| R_\psi (t) \|_{H^\s} & \leq C \e^2 \big(\| \pa_t (\phi_\ell \psi_\ell ) \|_{H^\s} +  \| \Delta \psi_2 \|_{H^\s}     +   \| |\psi_\ell|^2 \psi_\ell \|_{H^\s} +   \| \psi_2 \phi_\ell \|_{H^\s} +  \e^2\| \psi_2 |\psi_\ell|^2 \|_{H^\s} \big) \\
                                & \leq C( \| (\phi_0, \phi_1) \|_{H^{\s+4}},  \| \psi_0 \|_{H^{\s+2}} ) \e^2.
    \end{split}
\end{align}
These estimates \eqref{est:Rphi-2} and \eqref{est:Rpsi-2} confirm that the residuals are of order $O(\e^2)$.

By the initial data choice \eqref{eq:e2} and \eqref{init:zero},
all $O(\e)$ terms in the initial error are eliminated, and the initial errors are reduced to
\begin{align*}
    r_{\Phi}(0,x) & = \Phi(0,x) - \Phi^{a} (0,x) =  - \e^{2}\begin{pmatrix}
                                                                0                         \\
                                                                {\mu} \pa_t \phi_{{\ell}} \\
                                                                \frac{ \lambda}{\mu} |\psi_{{\ell}} |^2
                                                            \end{pmatrix}(0, x), \\
    r_{\psi}(0,x) & = \psi(0,x) - \psi^{a} (0,x) = - \e^2 \psi_2(0,x).
\end{align*}
Applying the regularity estimates for the WKB components, we obtain
\begin{align*}
    \| r_\Phi \|_{H^\s} + \| r_\psi \|_{H^\s} \leq  C(\| \phi_0, \phi_1 \|_{H^{\s+4}}, \| \psi_0 \|_{H^{\s+2}})  \e^2.
\end{align*}
This completes the proof of Proposition \ref{prop:wkb-N2}.

\subsection{High regularity case: proof of Proposition \ref{prop:wkb-N4}} \label{prof-prop2}
To capture higher-order nonlinear interactions, we introduce the second-order profiles $(\phi_c, \psi_c)$ together with their algebraic terms and define the approximate solution $(\Phi^a, \psi^a)$ by the full expansions \eqref{wkb-Phi-final} and \eqref{wkb-psi-final}.

We first give the estimates for the components $\phi_\ell$, $\psi_\ell$, $\phi_c$, and $\psi_c$. Under initial data satisfying \eqref{init-norm-4}, the estimates for $\phi_\ell$ and $\psi_\ell$ follow directly from \eqref{eq:sobolev-unified} and \eqref{est-dispersive-unified}:
\ba\label{est:psil-4}
     & \| \phi_\ell(t) \|_{H^{\s+8}} \leq C \| (\phi_0, \phi_1) \|_{H^{\s+8}}, \quad
    \| \psi_\ell(t) \|_{H^{\s+4}} \leq C \| \psi_0 \|_{H^{\s+4}},                                                       \\
     & \| \phi_\ell(t) \|_{W^{\s+4,\infty}} \leq C (1+t)^{-\frac{d}{2}} \| (\phi_0, \phi_1) \|_{W^{\s+4,1}\cap H^{\s+4}}, \\
     & \| \psi_\ell(t) \|_{W^{\s,\infty}} \leq C (1+t)^{-\frac{d}{2}} \| \psi_0 \|_{W^{\s,1}\cap H^\s}. 
\ea
For the component $\phi_c(t)$, \eqref{formula-phic} together with the estimates for $\phi_\ell(t)$ gives
\ba\label{refined-phi}
     & \| \phi_c(t) \|_{H^{\s+4}} \leq C t \| \Delta^2 \phi_\ell(t) \|_{H^{\s+4}} \leq C (1+t) \| (\phi_0, \phi_1) \|_{H^{\s+8}},                                                                                    \\
     & \| \phi_c(t) \|_{W^{\s,\infty}} \leq C t \| \Delta^2 \phi_\ell(t) \|_{W^{\s,\infty}} \leq C t (1+t)^{-\frac{d}{2}} \| (\phi_0, \phi_1) \|_{W^{\s+4,1} \cap H^{\s+8}} \leq C (1+t)^{-\frac{d}{2}+1}.  
\ea
For the component $\psi_c(t)$, by the expression \eqref{formula-psic} and using the algebraic $H^\s$ norm:
\begin{align*} 
    \| \psi_c(t) \|_{H^{\s+4}} & \leq C \int_0^t \| |\psi_\ell|^2 \psi_\ell(s) \|_{H^{\s+4}} \D{s} \leq C \int_{0}^t \| \psi_\ell(s) \|_{H^{\s+4}}^3 \D{s} \leq C (1+t) \| \psi_0 \|_{H^{\s+4}}^3.
\end{align*}
Another estimate for $\psi_c(t)$ is given by
\begin{align*}
    \norm{ \psi_{c}(t) }_{H^\s} \leq C \int_0^t \norm{ |\psi_{\ell} |^2 \psi_{\ell}(s) }_{H^\s} \D{s} \leq C \int_0^t \norm{ \psi_{\ell}(s) }_{W^{\s,\infty}}^2 \norm{ \psi_{\ell}(s) }_{H^\s} \D{s}.
\end{align*}
Applying the decay estimate for $\psi_{\ell}(t)$ from \eqref{est:psil-4}, we obtain
\begin{equation} \label{refined-psi}
    \norm{ \psi_{c}(t) }_{H^\s} \leq C \norm{ \psi_0}_{W^{\s,1} \cap H^\s}^3 \int_0^t (1+s)^{-d}\D{s} \leq C (1 +t)^{-d+1}.
\end{equation}
This result shows a decay estimate but requires $W^{\s,1}$ regularity.

Now we give the estimates for all the WKB profiles, which consist of $\phi_\ell$, $\psi_\ell$, $\phi_c$, and $\psi_c$. The following estimates hold for all $t>0$:
\begin{align*}
    \| \Phi_{0,1}(t)\|_{H^{\s+8}}  & \leq C \| \phi_\ell(t) \|_{H^{\s+8}} \leq C \| (\phi_0, \phi_1) \|_{H^{\s+8}},                                            \\
    \| \Phi_{1,1}(t) \|_{H^{\s+7}} & \leq C \| \nabla \phi_\ell(t) \|_{H^{\s+7}} \leq C \| (\phi_0,\phi_1) \|_{H^{\s+8}},                                      \\
    \| \Phi_{2,0}(t) \|_{H^{\s+4}} & \leq  C \| |\psi_\ell(t)|^2 \|_{H^{\s+4}} \leq C \| \psi_0 \|_{H^{\s+4}}^2 ,                                              \\
    \| \Phi_{2,1}(t) \|_{H^{\s+4}} & \leq C \| \pa_t \phi_\ell(t) \|_{H^{\s+4}} + C \| \phi_c(t) \|_{H^{\s+4}} \leq C(1+t)  \| (\phi_0, \phi_1)\|_{H^{\s+8}} .
\end{align*}
Similarly, the remaining high-order components satisfy
\begin{align*}
    \| \Phi_{3,0}(t) \|_{H^{\s+3}} & \leq C \| \nabla \phi_c(t) \|_{H^{\s+3}} \leq C t \| \nabla \Delta^2 \phi_\ell(t) \|_{H^{\s+3}} \leq  C(1+t) \| (\phi_0, \phi_1) \|_{H^{\s+8}},
\end{align*}
and
\begin{align}\label{Phi4}
    \begin{split}
        \| \Phi_{4,0}(t) \|_{H^{\s+2}} & \leq  C \| \pa_t |\psi_\ell|^2(t) \|_{H^{\s+2}} + C \| \psi_\ell(t) \|_{H^{\s+2}} \| \psi_c(t) \|_{H^{\s+2}}                          \\
                                       & \leq  C \| \psi_0 \|_{H^{\s +4}}^2 + C(1+t) \| \psi_0 \|_{H^{\s+2}}^4,                                                                \\
        \| \Phi_{4,1}(t) \|_{H^{\s+2}} & \leq C\| \pa_t \phi_c(t) \|_{H^{\s+2}} \leq C t\| \Delta^3 \phi_\ell(t) \|_{H^{\s+2}} \leq C (1+t) \| (\phi_0, \phi_1) \|_{H^{\s+8}}.
    \end{split}
\end{align}

For the Schr\"odinger part,
\begin{align*}
    \| \psi_{0,0}(t) \|_{H^{\s+4}}                                  & \leq C \| \psi_\ell(t) \|_{H^{\s+4}} \leq C \| \psi_0 \|_{H^{\s+4}},                                                               \\
    \| \psi_{2,0}(t) \|_{H^{\s+4}}                                  & \leq  C \| \psi_c(t) \|_{H^{\s+4}} \leq C(1+t) \|\psi_0\|_{H^{\s+4}},                                                              \\
    \| \psi_{2,1}(t)\|_{H^{\s+4}} + \| \psi_{2,-1}(t) \|_{H^{\s+4}} & \leq C \| \phi_\ell(t) \|_{H^{\s+4}} \| \psi_\ell(t) \|_{H^{\s+4}} \leq C\|(\phi_0, \phi_1) \|_{H^{\s+4}} \| \psi_0 \|_{H^{\s+4}},
\end{align*}
and
\begin{align}  \label{psi4}
    \begin{split}
         & \| \psi_{4,1}(t) \|_{H^{\s+4}} + \| \psi_{4,-1}(t) \|_{H^{\s+4}}                                                                                                                                             \\
         & \quad \leq C \| \phi_c \psi_\ell(t) \|_{H^{\s+4}} + C \| \phi_\ell \psi_c(t) \|_{H^{\s+4}}    + C \| \pa_t (\phi_\ell \psi_\ell)(t) \|_{H^{\s+4}} + C \| \Delta (\phi_\ell \psi_\ell)(t) \|_{H^{\s+4}}       \\
         & \quad  \leq C(1+t) \| (\phi_0, \phi_1) \|_{H^{\s+8}}  \| \psi_0 \|_{H^{\s +4}} + C(1+t) \|(\phi_0, \phi_1) \|_{H^{\s+4}} \| \psi_0 \|_{H^{\s+4}}     ,                                                       \\
         & \| \psi_{4,2}(t) \|_{H^{\s+4}} + \| \psi_{4,-2}(t) \|_{H^{\s+4}}   \leq C \| \phi_\ell(t) \|_{H^{\s+4}}^2 \|\psi_{\ell}(t) \|_{H^{\s+4}} \leq C \| (\phi_0, \phi_1) \|_{H^{\s+4}}^2 \| \psi_0 \|_{H^{\s+4}}.
    \end{split}
\end{align}

From the WKB expansion in Section \ref{sec:profiles},  it follows that the approximate solution $(\Phi^{a}, \psi^{a})$ defined in \eqref{wkb-Phi-final} and \eqref{wkb-psi-final} satisfies the the approximate system \eqref{app-sys} with the Klein-Gordon residual
\begin{equation*}
    \begin{split}
        R_{\Phi}  = \begin{pmatrix}
                        \frac{i\e^3 \l}{\mu} \nabla \pa_t |\psi_{\ell}|^2                                                   \\
                        \frac{2\e^4}{\mu} \Re( \pa_t^2 \phi_{c} e^{i\th} ) +  \frac{i \e^4 \l}{\mu} \pa_t^2 |\psi_{\ell}|^2 \\
                        \frac{2\e^4 \l}{\mu} \pa_t  \Re( \psi_{\ell} \overline{\psi_{c}})
                    \end{pmatrix}                                                                     - \lambda \e^4 \begin{pmatrix}
                                                                                                                         0                                                                                                             \\
                                                                                                                         |\psi_2|^2  +2\Re(\psi_4 \overline{\psi_{\ell}})   + 2 \e^2  \Re( \psi_4 \overline{\psi_2}) + \e^4 |\psi_4|^2 \\
                                                                                                                         0
                                                                                                                     \end{pmatrix},
    \end{split}
\end{equation*}
and the Schr\"odinger residual
\begin{equation*} 
    \begin{split}
        R_{\psi} & =   \e^4 (\pa_t \psi_{4,1}e^{i \th} + \pa_t \psi_{4,2} e^{2i\th} +\pa_t \psi_{4,-1} e^{-i \th} + \pa_t \psi_{4,-2} e^{-2i \th}  - i \Delta \psi_4)
        - \e^4 \frac{i \l}{\mu} \Big[ (\phi_{\ell} e^{i \th} + \overline{\phi_{\ell}} e^{-i\th}) \psi_4                                                              \\
                 & \quad + (  \phi_{c} e^{i \th} +\overline{\phi_{c}} e^{-i \th} +\frac{\l}{\mu} |\psi_{\ell}|^2 ) (\psi_2 +\e^2 \psi_4)
            + \frac{2\l}{\mu} \Re(\psi_{\ell} \overline{\psi_{c}})(\psi_{\ell} + \e^2 \psi_2 + \e^4 \psi_4 )   \Big].
    \end{split}
\end{equation*}
Under the initial assumption in \eqref{init-norm-4}, the residuals are bounded by
\begin{align*}
    \norm{ R_{\Phi}  (t) }_{H^\s(\R^d)} \leq C \e^3  ( 1+ \e t), \quad \norm{ R_{\psi}  (t) }_{H^\s(\R^d)} \leq C \e^4 (1+t).
\end{align*}
For the initial errors, substituting the explicit profiles into the definition yields
\begin{align*}
    r_{\Phi}(x) & = -\e^2 \begin{pmatrix}
                              0                                                                                  \\
                              \frac{2}{\mu} \Re( \pa_t \phi_{\ell} )(0,x) + \e^2
                              \frac{1}{\mu} \pa_t |\psi_{\ell}(0,x)|^2 + \frac{2}{\mu} \Re(\pa_t \phi_{c} )(0,x) \\
                              \frac{\l}{\mu} |\psi_{\ell}(0,x)|^2 + \e^2 \frac{2\l}{\mu} \Re(\psi_{\ell} \overline{\psi_{c}}) (0,x)
                          \end{pmatrix}, \\
    r_{\psi}(x) & =  -2 \e^2 \Re \left( \psi_{2,1}(0,x) + \psi_{2,-1}(0,x) \right) - 2\e^4 \Re(\psi_4(0,x) ) ,
\end{align*}
and these satisfy
\begin{equation*}
    \norm{ r_{\Phi}  }_{H^\s}  + \norm{ r_{\psi} }_{H^\s} \leq C \e^2,
\end{equation*}
where $C$ depends on $\| (\phi_0, \phi_1 )\|_{H^{\s+8}}$ and  $\| \psi_0 \|_{H^{\s+4}}$.

\begin{remark}
    In fact, the initial error can be reduced to $O(\e^3)$ by choosing non-zero initial data for $(\phi_c, \psi_c)$ to cancel the $O(\e^2)$ terms, as in previous works \cite{BLZ24, LZ25}.  However, this will not improve the final $O(\e^{2})$ convergence rates, as explained at the end of Section \ref{sec:reform}.
\end{remark}

\section{Stability and convergence rates} 

In this section, we establish the convergence rates by analyzing the stability of the WKB approximations.

\subsection{Perturbation system}
We introduce the perturbations
\begin{equation} \label{def-dotu}
    \dot{\Phi}(t,x) := \Phi(t,x) - \Phi^a(t,x), \quad \dot{\psi}(t,x): = \psi(t,x) - \psi^a(t,x), \quad \dot{u}(t,x): = (\dot{\Phi}(t,x), \dot{\psi}(t,x))^{\rm T}.
\end{equation}
Substituting \eqref{def-dotu} into \eqref{4-sys-kgs} and using \eqref{app-sys} yields the perturbation system
\begin{equation} \label{dotu-sys-0}
    \begin{dcases}
        \partial_t \dot{\Phi} + \frac{1}{\e} A(\partial_x) \dot{\Phi} + \frac{1}{\e^2} A_0 \dot{\Phi} = \mathcal{R}_\Phi, \\
        \partial_t \dot{\psi} - i \Delta \dot{\psi} = \mathcal{R}_\psi,                                                   \\
        \dot{\Phi}(0) = r_{\Phi}, \quad \dot{\psi}(0) = r_{\psi},
    \end{dcases}
\end{equation}
where the residual $(\mathcal{R}_\Phi, \mathcal{R}_\psi)$ consists of the nonlinear difference and the WKB residuals:
\begin{equation} \label{def:macR}
    \mathcal{R}_\Phi := F(\psi, \psi) - F(\psi^a, \psi^a) - R_\Phi, \quad  \mathcal{R}_\psi := g(\Phi, \psi) - g(\Phi^a, \psi^a) - R_\psi.
\end{equation}

Rewriting the system \eqref{dotu-sys-0} in terms of the perturbation $\dot{u}$ yields
\begin{align} \label{dotu-sys}
    \begin{dcases}
        \pa_t \dot{u} + \frac{i}{\e^2} \A(\e, D) \dot{u} = (\mathcal{R}_\Phi, \mathcal{R}_\psi)^{\rm T}, \\
        \dot{u}(0) = (r_\Phi, r_\psi)^{\rm T} ,
    \end{dcases}
\end{align}
where the self-adjoint operator $\mathcal{A}(\e, \nabla)$ is given by
\begin{equation*} 
    \mathcal{A}(\e, \nabla) := \begin{pmatrix}
        {0}_{d \times d}   & i\e \nabla & {0}_{d \times 1} & 0            \\
        i\e \nabla^{\rm T} & 0          & -i \mu           & 0            \\
        {0}_{1 \times d}   & i \mu      & 0                & 0            \\
        0                  & 0          & 0                & -\e^2 \Delta
    \end{pmatrix},
\end{equation*}
and consequently
\begin{equation} \label{eq:conversedHs}
    \norm{ e^{-\frac{i}{\e^2} \mathcal{A}(\e, \nabla) t}  \dot{u}(t) }_{H^\sigma} = \norm{ \dot{u}(t) }_{H^\sigma}, \quad \forall \, t >0.
\end{equation}
By standard theory for symmetric hyperbolic systems, for any $\s > \frac{d}{2}$,   there exists a unique solution $\dot{u} \in C([0, T_\e^*); H^\sigma(\mathbb{R}^d))$ to \eqref{dotu-sys}, where $T_\e^{*}$ is the existence time and will to be determined later.

Applying Duhamel's principle, the solution satisfies for all $t \in [0, T_\e^*)$:
\begin{equation*}
    \dot{u}(t) = e^{-\frac{i}{\e^2} \mathcal{A}(\e, \nabla) t} \dot{u}(0) + \int_0^t e^{-\frac{i}{\e^2} \mathcal{A}(\e, \nabla) (t-s)} \big( \mathcal{R}_\Phi, \mathcal{R}_\psi \big)^{\rm T} (s) \D{s}.
\end{equation*}
By \eqref{def:macR} and \eqref{eq:conversedHs}, we obtain
\begin{align} \label{energy-expan}
    \begin{split}
        \norm{ \dot{u}(t)}_{H^\s}
         & \leq  \| (r_\Phi, r_\psi) \|_{H^\s} + \int_0^t \norm{ (R_\Phi, R_\psi)(s) }_{H^\s} \D{s}                                                 \\
         & \quad + \int_0^t \norm{ F(\psi, \psi)(s) - F(\psi^a, \psi^a)(s) }_{H^\s} + \norm{ g(\Phi, \psi)(s) - g(\Phi^a, \psi^a)(s) }_{H^\s}\D{s}.
    \end{split}
\end{align}
Propositions \ref{prop:wkb-N2} and \ref{prop:wkb-N4} show that the $H^\sigma$ norms of $(r_\Phi, r_\psi)$ and $(R_\Phi, R_\psi)$ are small. The main task is to control the errors $F(\psi, \psi) - F(\psi^a, \psi^a)$ and $g(\Phi, \psi) - g(\Phi^a, \psi^a)$ in energy inequality \eqref{energy-expan}.

\subsection{Proof of Theorem \ref{main-thm-1}} 

Under the assumptions on the initial data in Theorem \ref{main-thm-1}, we can construct a WKB solution of the form \eqref{sche-2} satisfying the properties given in  Proposition \ref{prop:wkb-N2}. Consequently, we have
\begin{align} \label{res-2}
    \norm{ (r_\Phi, r_\psi) }_{H^\s}+ \int_0^t \norm{ (R_\Phi, R_\psi)(s) }_{H^\s} \D{s} \leq C(1+t) \e^2.
\end{align}
To complete the energy estimate via \eqref{energy-expan}, it remains to bound the error
\begin{align*}
    \int_0^t (\norm{ F(\psi, \psi)(s) - F(\psi^a, \psi^a)(s) }_{H^\s} + \norm{ g(\Phi, \psi)(s) - g(\Phi^a, \psi^a)(s) }_{H^\s}) \D{s}.
\end{align*}
Using the form of $F$ in \eqref{4-defF} and $\psi^a$ in \eqref{sche-2} gives
\begin{align} \label{est:F1.3}
    \begin{split}
        F(\psi, \psi) - F(\psi^a, \psi^a) & =  2 F(\psi^a, \dot{\psi})   + F(\dot{\psi}, \dot{\psi})                                                                                                                 \\
                                          & = 2 F(\psi_\ell, \dot{\psi}) + \frac{2\l \e^2}{\mu^2}  F( \phi_\ell \psi_\ell e^{i\th} -  \xwl{\phi}_\ell \psi_\ell e^{-i\th}, \dot{\psi}) +  F(\dot{\psi}, \dot{\psi}).
    \end{split}
\end{align}
Thus,
\begin{align*}
     & \|  F(\psi, \psi) - F(\psi^a, \psi^a) \|_{H^\s}    \leq C \| \psi_\ell  \dot{\psi} \|_{H^\s} + C \e^2 \| \phi_\ell \psi_\ell \dot{\psi}  \|_{H^\s}  + C \e^2 \| \overline \phi_\ell \psi_\ell \dot{\psi}  \|_{H^\s}   + C \|  |\dot{\psi}|^2 \|_{H^\s}.
\end{align*}
By the bilinear estimate \eqref{bilinear-est-app} and the algebraic $H^\s$ norm, we obtain
\begin{align} \label{est-fdiff}
    \begin{split}
         & \| F(\psi, \psi)(t) - F(\psi^a, \psi^a)(t) \|_{H^\s}                                                                                                                                  \\
         & \quad \leq C \| \psi_\ell \|_{W^{\s, \infty}} \| \dot{\psi} \|_{H^\s} + C \e^2 \| \phi_\ell \|_{H^\s} \| \psi_\ell \|_{H^\s} \| \dot{\psi}  \|_{H^\s} + C  \|  \dot{\psi} \|_{H^\s}^2 \\
         & \quad \leq C(\| (\phi_0, \phi_1, \psi_0) \|_{W^{\s,1} \cap H^\s}) ( (1+t)^{-\frac{d}{2}} + \e^2  + \|\dot{u}(t) \|_{H^\s} ) \| \dot{u}(t) \|_{H^\s}.
    \end{split}
\end{align}

Similarly, noting that the function $g$ in \eqref{4-defF} only depends on the third component $\Phi^{(3)}$, the difference is given by
\begin{align} \label{est:g1.3}
    \begin{split}
         & g(\Phi,\psi) - g(\Phi^a, \psi^a)  = g(\dot{\Phi}, \psi_0) + g(\Phi_0, \dot{\psi})  +
        \e^2 g(\dot{\Phi}, \psi_2)      + \e^2 g(\Phi_2 , \dot{\psi}) + g(\dot{\Phi}, \dot{\psi})                                                                                                                                                                                                                                                 \\
         & =        \frac{i\l}{\mu} \big( \dot{\Phi}^{(3)} \psi_{\ell} + (\phi_{\ell} e^{i\th} + \overline{\phi_{\ell}} e^{-i\th}) \dot{\psi} + \e^2( \dot{\Phi}^{(3)} (\frac{\l}{\mu} \phi_\ell \psi_\ell - \frac{\l}{\mu} \xwl{ \phi_\ell} \psi_\ell e^{-i\th})  + \frac{\l}{\mu} |\psi_\ell|^2 \dot{\psi}) + \dot{\Phi}^{(3)}\dot{\psi} \big).
    \end{split}
\end{align}
Therefore,
\begin{align} \label{est-gdiff}
    \begin{split}
         & \| g(\Phi,\psi)(t) - g(\Phi^a, \psi^a)(t) \|_{H^\s}                                                                                                                                                                                                              \\
         & \quad \leq C\big( \| \dot{\Phi}^{(3)} \psi_{\ell}\|_{H^\s} +   \| \phi_{\ell}  \dot{\psi} \|_{H^\s} +  \e^2 \|  \dot{\Phi}^{(3)} \phi_\ell \psi_\ell \|_{H^\s} +  \e^2 \| |\psi_\ell|^2 \dot{\psi}  \|_{H^\s} +   \| \dot{\Phi}^{(3)}\dot{\psi}  \|_{H^\s} \big) \\
         & \quad \leq C (\| (\phi_0, \phi_1, \psi_0) \|_{W^{\s,1} \cap H^\s}) ( (1+t)^{-\frac{d}{2}} + \e^2 + \|\dot{u}(t) \|_{H^\s} ) \| \dot{u}(t) \|_{H^\s}.
    \end{split}
\end{align}

Substituting \eqref{res-2}, \eqref{est-fdiff} and \eqref{est-gdiff} into \eqref{energy-expan} yields
\begin{equation} \label{energy-case2}
    \norm{ \dot{u}(t) }_{H^\s} \leq C(1+t) \e^2 + C \int_0^t ( (1+s)^{-\frac{d}{2}}  + \|\dot{u} (s)\|_{H^\s} ) \| \dot{u} (s)\|_{H^\s}  \D{s}.
\end{equation}

For the fixed constant $C$ in \eqref{energy-case2}, we employ a bootstrap argument to show the long-time estimates.
Define the large constant $M_{1}$ and small time $T_{1}$ as
\begin{align*} 
    M_1 := (C+1)\exp(\frac{C d}{d-2}),\quad T_{1} := \frac{1}{2M_{1}},
\end{align*}
and the bootstrap time $T_\e^1$:
\begin{equation*} 
    T_\e^1 := \sup \big\{  t \in [0, \min\{ T_\e^*, \frac{T_1}{\e}\} ) : \sup_{0 \leq s \leq t} \norm{\dot{u}(s) }_{H^\s} \leq M_1 (1+t)\e^2 \big\}.
\end{equation*}

Let $t <  T_\e^1$. Under the bootstrap hypothesis, \eqref{energy-case2} becomes
\begin{equation*}
    \norm{ \dot{u}(t) }_{H^\s} \leq C(1+t) \e^2 + C \int_0^t \big( (1+s)^{-\frac{d}{2}}  +  M_1 (1+s) \e^2 \big) \norm{ \dot{u}(s) }_{H^\s} \D{s}.
\end{equation*}
By Gr\"onwall's inequality, we obtain
\begin{align} \label{term:exp-2}
    \begin{split}
        \sup_{0\leq s \leq t}\norm{ \dot{u}(s) }_{H^\s} & \leq C(1+t) \e^2 \exp \big( C \int_0^t \big( (1+s)^{-\frac{d}{2}} +  M_1 (1+s) \e^2 \big) \D{s} \big) \\
                                                        & \leq C (1+t)\e^2 \exp(C (\frac{2}{d-2}  +  M_1 \e^2 (t + \frac{1}{2}t^2) ) .
    \end{split}
\end{align}
Recall $T_1 = \frac{1}{2M_1} <1$. Then, the quantity in \eqref{term:exp-2} is bounded by
\begin{equation*}
    \frac{2}{d-2}  +  M_1 \e^2 (t + \frac{1}{2}t^2)  \leq    \frac{2}{d-2} + M_1 ( \e T_1 + \frac{1}{2} T_1^2 ) \leq  \frac{2}{d-2} +  2 M_1 T_1 \leq    \frac{d}{d-2} .
\end{equation*}
Thus, for all $t < T_\e^1 $, we have
\begin{equation*}
    \sup_{0\leq s\leq t}\norm{ \dot{u}(s) }_{H^\s} \leq C(1+t) \e^2 \exp(\frac{C d}{d-2}) \leq \frac{C}{C+1} M_1 (1+t) \e^2.
\end{equation*}

A standard continuity argument implies that $T_\e^1$ cannot be smaller than $\min\{ \frac{T_1}{\e}, T_\e^*\}$. Furthermore, since the $H^\sigma$ norm of function $\dot{u}$ does not blow up as $t \to \frac{T_1}{\e}$, the local-in-time existence theory for \eqref{dotu-sys} ensures that the maximal existence time  $T_\e^* \ge \frac{T_1}{\e}$, which in turn implies $T_\e^1 \ge \frac{T_1}{\e}$. We thus conclude that
\begin{equation*} 
    \sup_{0 \leq t \leq \frac{T_1}{\e}} \norm{ \dot{u}(t) }_{H^\s} \leq M_1 (1+t) \e^2.
\end{equation*}

Returning to the components $\Phi$, $\Phi^a$, $\psi$ and $\psi^a$,
we obtain
\begin{align} \label{boot-res-2}
    \norm{\Phi(t) - \Phi^a(t)}_{H^\sigma} + \norm{\psi(t) - \psi^a(t)}_{H^\sigma} \leq C(1+t)\e^2,
\end{align}
which is the desired error estimate in \eqref{eq:u-ua-2}.

By the WKB construction \eqref{sche-2}, we find that
\begin{align*}
    (\Phi^a)^{(3)}(t) - ( e^{\frac{i\mu t}{\e^2}} \phi_\ell(t) + e^{-\frac{i\mu t}{\e^2} } \xwl{\phi_\ell(t)} ) & =\frac{\e^2 \l}{\mu} |\psi_\ell|^2 ,            \quad
    \psi^a(t) - \psi_\ell(t)     =  \e^2 \big(   \frac{\lambda}{\mu^2} \phi_{\ell} \psi_{\ell}  e^{i \theta} - \frac{\l}{\mu^2} \overline{\phi_{\ell}} \psi_{\ell} e^{- i\theta} \big).
\end{align*}
Thus,
\begin{align} \label{diff-2}
    \| (\Phi^{a})^{(3)}(t) - (  e^{\frac{i \mu t}{\e^2}} \phi_{\ell}(t)  +  e^{- \frac{i \mu t}{\e^2}} \overline{\phi_{\ell}}(t) ) \|_{H^\s} +  \| \psi^{a}(t) - \psi_{\ell}(t) \|_{H^\s} \leq C \e^2.
\end{align}
Therefore, by \eqref{boot-res-2} and \eqref{diff-2}, the triangle inequality gives
\begin{equation*}
    \| \mu \phi(t) - (  e^{\frac{i \mu t}{\e^2}} \phi_{\ell}(t)  +  e^{- \frac{i \mu t}{\e^2}} \overline{\phi_{\ell}}(t) ) \|_{H^\s} +  \| \psi(t) - \psi_{\ell}(t) \|_{H^\s} \leq C (1+t)\e^2 + C\e^2 \leq C (1+t)\e^2.
\end{equation*}
This completes the proof of Theorem \ref{main-thm-1}.

\subsection{Proof of Theorem \ref{main-thm-2}}

The proof is based on the energy inequality {energy-expan}. Under the regularity assumptions on the initial data in Theorem \ref{main-thm-2}, we can construct a WKB solution given by \eqref{wkb-Phi-final} and \eqref{wkb-psi-final} satisfying the properties given in  Proposition \ref{prop:wkb-N4}. Thus,
\begin{align} \label{res-4}
    \norm{ (r_\Phi, r_\psi) }_{H^\s}+ \int_0^t \norm{ (R_\Phi, R_\psi)(s) }_{H^\s} \D{s} \leq C \e^2 + \int_0^t C(\e^3 + \e^4 s) \D{s} \leq C(1+ \e t + \e^2 t^2) \e^2.
\end{align}

To capture the error $F(\psi, \psi) - F(\psi^a, \psi^a)$, we use the form of $F$ in \eqref{4-defF} and $\psi^a$ in \eqref{wkb-psi-final} to obtain
\begin{align*} 
     & F(\psi, \psi) - F(\psi^a, \psi^a)   = 2 F(\psi_0, \dot{\psi}) + 2\e^2 F(\psi_2 +\e^2 \psi_4, \dot{\psi}) +  F(\dot{\psi}, \dot{\psi})                                                                                                    \\
     & = 2 F(\psi_\ell, \dot{\psi}) + \frac{2\l \e^2}{\mu^2}  F( \phi_\ell \psi_\ell e^{i\th} -  \xwl{\phi}_\ell \psi_\ell e^{-i\th}, \dot{\psi}) + 2\e^2 F(\psi_c, \dot{\psi}) +   2\e^4 F(\psi_4, \dot{\psi}) +    F(\dot{\psi}, \dot{\psi}).
\end{align*}
Comparing with the previous case \eqref{est:F1.3}, the difference contains additional terms involving $\psi_c$ and $\psi_4$. These two terms are estimated as follows:
\begin{align*}
    \|  \e^2 F(\psi_c, \dot{\psi}) + \e^4 F(\psi_4, \dot{\psi})  \|_{H^\s}
     & \leq C \e^2 \| \psi_c \|_{H^\s}  \| \dot{\psi} \|_{H^\s} + C \e^4 \| \psi_4 \|_{H^\s} \| \dot{\psi} \|_{H^\s} .
\end{align*}

By the estimates for $\psi_4$ in \eqref{psi4} and $\psi_c$ in \eqref{refined-psi}, we have
\begin{align*}
    \| \e^2 F(\psi_c, \dot{\psi}) (t)+ \e^4 F(\psi_4, \dot{\psi})(t) \|_{H^{\s}} \leq C(\|(\phi_0, \phi_1 )\|_{H^{\s+8}}, \| \psi_0 \|_{W^{\s,1} \cap H^{\s+4}}) ((1+t)^{-d+1} \e^2 + t \e^4).
\end{align*}
Therefore,
\begin{align} \label{est:F4}
    \begin{split}
        \| F(\psi, \psi)(t) - F(\psi^a, \psi^a)(t) \|_{H^\s} & \leq C( (1+t)^{-\frac{d}{2}} + \e^2 +t \e^4 +  \|\dot{u}(t) \|_{H^\s} ) \| \dot{u} (t)\|_{H^\s}.
    \end{split}
\end{align}
Similarly, for the function $g$ defined in \eqref{4-defF}, we have
\begin{align*}
     & g(\Phi, \psi) - g(\Phi^a, \psi^a)                                                                                                                                                                                                                                                                                                           \\
     & =    \frac{i\l}{\mu} \big( \dot{\Phi}^{(3)} \psi_{\ell} + (\phi_{\ell} e^{i\th} + \overline{\phi_{\ell}} e^{-i\th}) \dot{\psi} + \e^2 \big( \dot{\Phi}^{(3)} (\frac{\l}{\mu} \phi_\ell \psi_\ell - \frac{\l}{\mu} \xwl{ \phi_\ell} \psi_\ell e^{-i\th})  + \frac{\l}{\mu} |\psi_\ell|^2 \dot{\psi} \big) + \dot{\Phi}^{(3)}\dot{\psi} \big) \\
     & \quad  + \e^2 \frac{i\l}{\mu} \dot{\Phi}^{(3)} \psi_c + \e^2 \frac{i\l}{\mu} (\phi_c e^{i\th} + \xwl{\phi_c} e^{-i\th}) \dot{\psi} + \e^4 \frac{i\l}{\mu} \big( \dot{\Phi}^{(3)} \psi_4 +   \Phi_4^{(3)}  \dot{\psi} \big).
\end{align*}
Compared to \eqref{est:g1.3},
only the last line contains additional terms. Thus,
\begin{align*}
     & \| \e^2 \frac{i\l}{\mu} \dot{\Phi}^{(3)} \psi_c + \e^2 \frac{i\l}{\mu} (\phi_c e^{i\th} + \xwl{\phi_c} e^{-i\th}) \dot{\psi} + \e^4 \frac{i\l}{\mu} \big( \dot{\Phi}^{(3)} \psi_4 +   \Phi_4^{(3)}  \dot{\psi} \big) \|_{H^\s} \\
     & \quad  \leq C \e^2 (\| \dot{\Phi}^{(3)} \psi_c \|_{H^\s} +  \|\phi_c \dot{\psi} \|_{H^\s} +  \|\xwl{\phi_c} \dot{\psi} \|_{H^\s}) + C \e^4 (\|   \dot{\Phi}^{(3)}\psi_4 \|_{H^\s} + \|   \Phi_4^{(3)} \dot{\psi} \|_{H^\s}).
\end{align*}
Moreover, by \eqref{refined-psi}, we have
\begin{align} \label{L2-4-1}
    \| \dot{\Phi}^{(3)} \psi_c(t) \|_{H^\s} \leq C  \| \dot{\Phi}(t) \|_{H^\s} \leq C \| \dot{u}(t) \|_{H^\s}.
\end{align}
Applying the bilinear estimate \eqref{bilinear-est-app} to $\phi_c \dot{\psi}$ and $\xwl{\phi_c} \dot{\psi}$ gives
\begin{equation} \label{L2-4-2}
    \norm{\phi_{c} \dot{\psi}(t) }_{H^\s} + \|\xwl{\phi_c} \dot{\psi}(t) \|_{H^\s}  \leq C  \norm{ \phi_c(t)  }_{W^{\s, \infty}} \norm{\dot{\psi}(t) }_{H^\s}.
\end{equation}
By the estimate for $\phi_c$ in \eqref{refined-phi}, we have
\begin{align*}
    \norm{\phi_{c} \dot{\psi}(t) }_{H^\s} + \|\xwl{\phi_c} \dot{\psi}(t) \|_{H^\s} & \leq C(\norm{(\phi_{0}, \phi_1) }_{W^{\s+4,1} \cap H^{\s+4}} ) (1+t)^{-\frac{d}{2}+1} \norm{\dot{u}(t) }_{H^\s}.
\end{align*}

For the last terms involving $\Phi_4^{(3)}$ and $\psi_4$, we use \eqref{Phi4} and \eqref{psi4} to obtain
\begin{align}  \label{L2-4-3}
    \|  \psi_4 \dot{\Phi}^{(3)}(t) \|_{H^\s} + \| \Phi_4^{(3)} \dot{\psi}(t) \|_{H^\s} & \leq C(1+t) \| \dot{u}(t) \|_{H^\s}.
\end{align}
Consequently, combining \eqref{L2-4-1}, \eqref{L2-4-2}, \eqref{L2-4-3} and the estimate \eqref{est-gdiff}, it follows that
\begin{align} \label{est:g4}
    \begin{split}
         & \norm{ g(\Phi, \psi)(t) - g(\Phi^a,\psi^a)(t) }_{H^\s}                                                                                                \\
         & \quad  \leq  C \big( (1+t)^{-\frac{d}{2}} + \e^2 + \|\dot{u}(t) \|_{H^\s}  + \e^2 + \e^2(1+t)^{-\frac{d}{2}+1} + \e^4 (1+t) ) \| \dot{u}(t) \|_{H^\s} \\
         & \quad \leq C ( (1+t)^{-\frac{d}{2}} + \e^2(1+\e^2 t) + \| \dot{u}(t) \|_{H^\s} \big) \| \dot{u}(t) \|_{H^\s}.
    \end{split}
\end{align}
Taking the estimates \eqref{res-4}, \eqref{est:F4} and \eqref{est:g4} into the energy inequality \eqref{energy-expan}, we obtain
\begin{equation*} 
    \norm{ \dot{u}(t) }_{H^\s} \leq C(1 + \e t + \e^2 t^2)\e^2  + C \int_0^t \big( {(1+s)^{-\frac{d}{2}}} + \e^2(1+\e^2 s) + \| \dot{u}(s) \|_{H^\s} \big) \|\dot{u}(s) \|_{H^\s} \D{s},
\end{equation*}
which holds for all $t$ in the maximal existence time interval $[0, T_\e^*)$.

Define the large constant $M_{2}$ and small time $T_{2}$ as
\begin{align*} 
    M_2 :=   (C+1) \exp( \frac{C d}{d-2} ), \quad T_2 := \frac{1}{2+3 M_2},
\end{align*}
and the bootstarp time
\begin{equation*} 
    T_\e^2 := \sup \big\{ t \in [0,  \min\{ T_\e^*, \frac{T_2}{\e^{\frac{4}{3}}} \}) : \sup_{0 \leq s \leq t} \norm{ \dot{u}(s) }_{H^\s} \leq M_2 (1 + \e t + \e^2 t^2) \e^2 \big\}.
\end{equation*}

Let $t< T_\e^2$. Under the bootstrap hypothesis,
the estimate for $\dot{u}(t)$ becomes
\begin{align} \label{eq:energy-2}
    \begin{split}
        \sup_{0\leq s\leq t} \| \dot{u}(s) \|_{H^\s} & \leq C(1 + \e t + \e^2 t^2)\e^2                                                                                                         \\
                                                     & \quad   + \int_0^t \big( {(1+s)^{-\frac{d}{2}}} + \e^2(1+\e^2 s) +  M_2 (1 + \e s + \e^2 s^2)\e^2 \big) \| \dot{u}(s) \|_{H^\s}  \D{s}.
    \end{split}
\end{align}
Applying Gr\"onwall's inequality to \eqref{eq:energy-2} gives
\begin{align*}
    \sup_{0\leq s\leq t}  \norm{ \dot{u}(s) }_{H^\s}  \leq   C(1 + \e t + \e^2 t^2)\e^2   \exp (  C \int_0^t ( {(1+s)^{-\frac{d}{2}}} + \e^2(1+\e^2 s) + M_2(1 + \e s + \e^2 s^2)\e^2 ) \D{s} ).
\end{align*}

This implies for all $t  < T_\e^2 $ that
\begin{align*} 
    \begin{split}
        \sup_{0\leq s\leq t}  \| \dot{u}(s) \|_{H^\s} & \leq C(1 + \e t + \e^2 t^2)\e^2   \exp\big( C(\frac{2}{d-2} + T_2 \e^{\frac{2}{3}} + \f{2} T_2 \e^{\frac{4}{3}} + M_2  ( T_2 \e^{\frac{2}{3}}   + \f{2} \e^{\frac{1}{3} } T_2^2 +\f{3} T_2^3 ) ) \big).
    \end{split}
\end{align*}
Recalling $T_2 = \frac{1}{2 + 3 M_2}$, we thus have
\begin{align*}
     & T_2 \e^{\frac{2}{3}} + \f{2} T_2 \e^{\frac{4}{3}} + M_2  ( T_2 \e^{\frac{2}{3}}   + \f{2} \e^{\frac{1}{3} } T_2^2 +\f{3} T_2^3 ) \leq 2 T_2 + 3 M_2 T_2 \leq 1.
\end{align*}
Hence, the integral in the exponent is uniformly bounded by $\frac{C d}{d-2}$. Consequently,
\begin{align*}
    \sup_{0\leq s\leq t} \| \dot{u}(s) \|_{H^\s} \leq C(1+\e t+ \e^2 t^2) \e^2 \exp( \frac{C d}{d-2} ) = \frac{C}{C+1} M_2 (1+\e t+\e^2 t^2) \e^2.
\end{align*}

By continuity and blow-up criterion for symmetric hyperbolic system, we conclude that $T_\e^* \ge  \frac{T_2}{\e^{\frac{4}{3}}} $ and $T_\e^2 \ge \frac{T_2}{\e^{\frac{4}{3}}}$.
Thus, the bootstrap estimate holds on the extended time interval  $t \in [0, \frac{T_2}{\e^{\frac{4}{3}}}]$:
\begin{equation*}
    \norm{ \dot{u}(t) }_{H^\s} \leq M_2 (1 + \e t + \e^2 t^2) \e^2.
\end{equation*}
Returning to the original variables, for every $t \in [0, \frac{T_2}{\e^{\frac{4}{3}}}]$, we obtain
\begin{equation} \label{eq:dotu4}
    \norm{ \Phi(t) - \Phi^a(t) }_{H^\s} + \norm{ \psi(t) - \psi^a(t) }_{H^\s} \leq C(1 + \e t + \e^2 t^2)\e^2.
\end{equation}

The difference between
$((\Phi^a)^{(3)} , \psi^a(t) )$ and $( e^{\frac{i\mu t}{\e^2}} \phi_\ell(t) + e^{-\frac{i\mu t}{\e^2} } \xwl{\phi_\ell(t)} , \psi_\ell(t))$ is derived from the construction \eqref{wkb-Phi-final} and \eqref{wkb-psi-final}:
\ba
& (\Phi^a)^{(3)}(t) -(  e^{\frac{i \mu t}{\e^2}} \phi_{\ell}(t)  +  e^{- \frac{i \mu t}{\e^2}} \overline{\phi_{\ell}}(t) ) = \frac{\e^2 \l}{\mu} |\psi_\ell|^2 + \e^2 ( \phi_c e^{i\th} + \xwl{\phi_c} e^{-i\th} ) + \e^4 \frac{2\l}{\mu} \Re(\psi_\ell \xwl{\psi_c}),\\
&    \psi^a(t) - \psi_\ell(t)  = \e^2 \big( \psi_c + \frac{\lambda}{\mu^2} \phi_{\ell} \psi_{\ell}  e^{i \theta} - \frac{\l}{\mu^2} \overline{\phi_{\ell}} \psi_{\ell} e^{- i\theta} \big) +\e^4 \psi_4.
\nn\ea
From the WKB component estimates in Proposition \ref{prop:wkb-N4}, we directly obtain
\ba
\norm{(\Phi^a)^{(3)}(t) -(  e^{\frac{i \mu t}{\e^2}} \phi_{\ell}(t)  +  e^{- \frac{i \mu t}{\e^2}} \overline{\phi_{\ell}}(t) ) }_{H^\s} \leq  C \e^2 (1 + t) + C \e^4 (1 +  t) \leq C \e^2 (1+t +\e^2 + \e^2 t).
\nn \ea
Combining this with \eqref{eq:dotu4}, we obtain
\begin{equation} \label{est:final1}
    \norm{ \mu \phi(t) - (  e^{\frac{i \mu t}{\e^2}} \phi_{\ell}(t)  +  e^{- \frac{i \mu t}{\e^2}} \overline{\phi_{\ell}}(t) ) }_{H^\s} \leq C(1 + t + \e^2 t^2) \e^2.
\end{equation}

For the Schr\"odinger component $\psi$, using \eqref{refined-psi} gives
\begin{equation} \label{eq:psia-psil}
    \| \psi^a(t) - \psi_\ell(t) \|_{H^\s} \leq C \e^2 + C \e^4 (1+t ) \leq C (1 + \e^2 +\e^2 t)  \e^2.
\end{equation}
Combining \eqref{eq:psia-psil} with the bootstrap estimate \eqref{eq:dotu4}, we finally obtain
\begin{equation} \label{est:final2}
    \norm{ \psi(t) - \psi_\ell(t) }_{H^\s} \leq C(1 + \e t + \e^2 t^2) \e^2.
\end{equation}
The above estimates \eqref{est:final1} and \eqref{est:final2} yield the desired error estimates  in Theorem \ref{main-thm-2}.

\section*{Acknowledgements}
W. B. is partially supported by the Ministry of Education of Singapore under its AcRF Tier 1 funding grant A-8003584-00-00.  Y. L. and Z. Z. are partially supported by Basic Research Program of Jiangsu (Grant No. BK20240058), and Y. L. is also partially supported by the Fundamental Research Funds for the Central Universities (Grant No. 2026300394).


\begin{thebibliography}{99}

    \bibitem{BC77}
    J.-B. Baillon and J. M. Chadam,
    The Cauchy problem for the coupled Schr\"odinger-Klein-Gordon equations,
    in \textit{Contemporary developments in continuum mechanics and partial differential equations},
    North-Holland Math. Stud., \textbf{30}, pp.~37--44,

    \bibitem{BFV15}
    C. Banquet, L. C. F. Ferreira, and E. J. Villamizar-Roa,
    On existence and scattering theory for the Klein-Gordon-Schr\"odinger system in an infinite $L^2$-norm setting,
    Ann. Mat. Pura Appl. (4) \textbf{194} (2015), no.~3, 781--804.

    \bibitem{BCZ14}
    W. Bao, Y. Cai and X. Zhao,
    A uniformly accurate multiscale time integrator pseudospectral method for the Klein-Gordon equation in the nonrelativistic limit regime, SIAM J. Numer. Anal. \textbf{52} (2014), 2488-2511.


    \bibitem{BL26}
    W. Bao and C. Liu,
    A uniformly accurate multiscale time integrator for the
    nonlinear Klein-Gordon equation in the nonrelativistic regime via simplified transmission conditions,
    arXiv: 2602.04988.

    \bibitem{BLZ24}
    W. Bao, Y. Lu, and Z. Zhang,
    Convergence rates in the nonrelativistic limit of the cubic Klein-Gordon equation,
    SIAM J. Math. Anal. \textbf{56} (2024), no.~5, 6822--6860.


    \bibitem{BZ17}
    W. Bao and X. Zhao,
    A uniformly accurate multiscale time integrator Fourier pseudospectral method for the Klein-Gordon-Schr\"odinger equations in the nonrelativistic limit regime,
    Numer. Math. \textbf{135} (2017), no.~3, 833--873.

    \bibitem{BZ19}
    W. Bao and X. Zhao,
    Comparison of numerical methods for the nonlinear Klein-Gordon equation in the nonrelativistic limit regime.
    J. Comput. Phys. \textbf{398} (2019), 108886, 30 pp.

    \bibitem{BMS04}
    P. Bechouche, N. J. Mauser, and S. Selberg,
    Nonrelativistic limit of Klein-Gordon-Maxwell to Schr\"odinger-Poisson,
    Amer. J. Math. \textbf{126} (2004), no.~1, 31--64.

    \bibitem{B99}
    J. Bourgain,
    \textit{Global solutions of nonlinear Schr\"odinger equations},
    Amer. Math. Soc. Colloq. Publ., vol.~46,
    American Mathematical Society, Providence, RI, 1999. viii+182 pp.


    \bibitem{CO06}
    Y. Cho and T. Ozawa,
    On the semirelativistic Hartree-type equation,
    SIAM J. Math. Anal. \textbf{38} (2006), no.~4, 1060--1074.

    \bibitem{FSZ24}
    C. Fan, G. Staffilani, and Z. Zhao,
    On decaying properties of nonlinear Schr\"odinger equations,
    SIAM J. Math. Anal., \textbf{56} (2024), pp. 3082--3109.

    \bibitem{FL26}
    Y. Feng and C. Liu, A uniformly accurate multiscale time integrator for the
    Klein-Gordon-Sch\"{o}dinger equations in the nonrelativistic regime via simplified transmission conditions,
    arXiv: 2605.12936.



    \bibitem{FT75}
    I. Fukuda and M. Tsutsumi,
    On coupled Klein-Gordon-Schr\"odinger equations. I,
    Bull. Sci. Engrg. Res. Lab. Waseda Univ. \textbf{69} (1975), 51--62.



    \bibitem{FT78}
    I. Fukuda and M. Tsutsumi,
    On coupled Klein-Gordon-Schr\"odinger equations. II,
    J. Math. Anal. Appl. \textbf{66} (1978), no.~2, 358--378.




    \bibitem{GO14}
    L. Grafakos and S. Oh,
    The Kato-Ponce inequality,
    Comm. Partial Differential Equations \textbf{39} (2014), no.~6, 1128--1157.

    \bibitem{GM95}
    B. Guo and C. Miao,
    Global existence and asymptotic behavior of solutions for the coupled Klein-Gordon-Schr\"odinger equations,
    Sci. China Ser. A. \textbf{38} (1995), no.~12, 1444--1456.

    \bibitem{HNR18}
    D. Han-Kwan, T. Nguyen and F. Rousset,
    Long time estimates for the Vlasov-Maxwell system in the non-relativistic limit,
    Comm. Math. Phys. \textbf{363} (2018), no.~2, 389--434.


    \bibitem{JMR00}
    J.-L. Joly, G. Métivier and J. Rauch,
    Transparent nonlinear geometric optics and Maxwell-Bloch equations,
    J. Differential Equations \textbf{166} (2000), no.~1, 175--250.


    \bibitem{KP88}
    T. Kato and G. Ponce,
    Commutator estimates and the Euler and Navier-Stokes equations,
    Comm. Pure Appl. Math. \textbf{41} (1988), no.~7, 891--907.


    \bibitem{LW25}
    Z. Lei and Y. Wu,
    Non-relativistic limit for the cubic nonlinear Klein-Gordon equations,
    arXiv: 2309.10235.




    \bibitem{LZ16}
    Y. Lu and Z. Zhang,
    Partially strong transparency conditions and a singular localization method in geometric optics,
    Arch. Ration. Mech. Anal. \textbf{222} (2016), no.~1, 245--283.

    \bibitem{LZ17}
    Y. Lu and Z. Zhang,
    Higher order asymptotic analysis of the Klein-Gordon equation in the non-relativistic limit regime,
    Asymptot. Anal. \textbf{102} (2017), no.~3-4, 157--175.



    \bibitem{LZ25}
    Y. Lu and Z. Zheng,
    Optimal convergence rates of the Klein-Gordon-Zakharov system in the non-relativistic limit,
    J. Differential Equations \textbf{431} (2025), Paper No. 113195, 52 pp.




    \bibitem{M01}
    S. Machihara,
    The nonrelativistic limit of the nonlinear Klein-Gordon equation,
    Funkcial. Ekvac. \textbf{44} (2001), no.~2, 243--252.


    \bibitem{M84}
    A. Majda,
    \textit{Compressible fluid flow and systems of conservation laws in several space variables},
    Appl. Math. Sci., vol.~53,
    Springer-Verlag, New York, 1984. viii+159 pp.


    \bibitem{MN02}
    N. Masmoudi and K. Nakanishi,
    From nonlinear Klein-Gordon equation to a system of coupled nonlinear Schr\"odinger equations,
    Math. Ann. \textbf{324} (2002), no.~2, 359--389.


    \bibitem{MN03}
    N. Masmoudi and K. Nakanishi,
    Nonrelativistic limit from Maxwell-Klein-Gordon and Maxwell-Dirac to Poisson-Schr\"odinger,
    Int. Math. Res. Not. \textbf{2003} (2003), no.~13, 697--734.

    \bibitem{MN08}
    N. Masmoudi and K. Nakanishi,
    Energy convergence for singular limits of Zakharov type systems,
    Invent. Math. \textbf{172} (2008), no.~3, 535--583.

    \bibitem{MN10}
    N. Masmoudi and K. Nakanishi,
    From the Klein-Gordon-Zakharov system to a singular nonlinear Schr\"odinger system,
    Ann. Inst. H. Poincar\'e C Anal. Non Lin\'eaire \textbf{27} (2010), no.~4, 1073--1096.


    \bibitem{MNO02}
    S. Machihara, K. Nakanishi, and T. Ozawa,
    Nonrelativistic limit in the energy space for nonlinear Klein-Gordon equations,
    Math. Ann. \textbf{322} (2002), no.~3, 603--621.



    \bibitem{M08}
    G. M\'etivier,
    \textit{Para-differential calculus and applications to the Cauchy problem for nonlinear systems},
    CRM Series, vol.~5,
    Edizioni della Normale, Pisa, 2008. xii+140 pp.


    \bibitem{MZ15}
    S. Missaoui and E. Zahrouni,
    Regularity of the attractor for a coupled Klein-Gordon-Schr\"odinger system with cubic nonlinearities in $\mathbb{R}^2$,
    Commun. Pure Appl. Anal. \textbf{14} (2015), no.~2, 695--716.


    \bibitem{P04}
    H. Pecher,
    Global solutions of the Klein-Gordon-Schr\"odinger system with rough data,
    Differential Integral Equations \textbf{17} (2004), no.~1-2, 179--214.

    \bibitem{P18}
    S. Pasquali,
    Almost global existence for the nonlinear Klein-Gordon equation in the nonrelativistic limit,
    J. Math. Phys. \textbf{59} (2018), no.~1, 011502, 15 pp.


    \bibitem{P19}
    S. Pasquali,
    Dynamics of the nonlinear Klein-Gordon equation in the nonrelativistic limit,
    Ann. Mat. Pura Appl. \textbf{4} 198 (2019), no.~3, 903--972.



    \bibitem{R12}
    J. Rauch,
    \textit{Hyperbolic partial differential equations and geometric optics},
    Grad. Stud. Math., vol.~133,
    American Mathematical Society, Providence, RI, 2012. xx+363 pp.


    \bibitem{SZ20}
    K. Schratz and X. Zhao,
    On comparison of asymptotic expansion techniques for nonlinear Klein-Gordon equation in the nonrelativistic limit regime,
    Discrete Contin. Dyn. Syst. Ser. B \textbf{25} (2020), no.~8, 2841--2865.


    \bibitem{T84}
    M. Tsutsumi,
    Nonrelativistic approximation of nonlinear Klein-Gordon equations in two space dimensions,
    Nonlinear Anal. \textbf{8} (1984), no.~6, 637--643.



    \bibitem{Y55}
    H. Yukawa,
    On the interaction of elementary particles. I,
    Progr. Theoret. Phys. Suppl. \textbf{1} (1955), 1--10.


\end{thebibliography}
\end{document}